\documentclass[10pt]{article}

\usepackage[dvips]{graphicx}
\usepackage{amsmath,amssymb,epsfig,subfigure}

\title{Some extensions of the class of convex bodies}
\author{
 Vladimir Golubyatnikov
 \footnote{The research was partially supported by leading scientific schools grant
 No. 8526.2006.1 of President of Russian Federation and by Fulbright visiting scholar program.}\\
 {\em\small Sobolev Institute of Mathematics, Novosibirsk, 630090, Russia} \\
 {\em\small e-mail: glbtn@math.\,nsc.\,ru}\\
 {\small and}\\
 Vladimir Rovenski
 \footnote{The research was supported by the Caesarea Edmond Benjamin de Rothschild Foundation Institute for Interdisciplinary Applications of Computer Science at the University of Haifa.}\\
 {\em\small
 University of Haifa,
 Haifa, 31905, Israel}\\
 {\em\small e-mail: rovenski@math.\,haifa.\,ac.\,il}
 }

\date{}

\begin{document}

\newtheorem{theo}{Theorem}
\newtheorem{defi}{Definition}
\newtheorem{lem}{Lemma}
\newtheorem{prob}{Problem}
\newtheorem{prop}{Proposition}
\newtheorem{cor}{Corollary}
\newtheorem{exam}{Example}
\newtheorem{rem}{Remark}
\newtheorem{exer}{Exercise}
\def\bbr{\mathbb{E}}
\def\bbn{\mathbb{N}}
\def\bbc{\mathbb{C}}
\def\bbh{\mathbb{H}}
\def\bbk{\mathbb{K}}
\def\bbz{\mathbb{Z}}
\def\bn{{\bf n}}
\def\ba{{\bf a}}
\def\br{{\bf r}}
\newcommand{\eps}{\varepsilon}

\maketitle

\begin{abstract} We introduce and study a new class of $\eps$-convex bodies (extending the class of convex bodies) in metric and normed linear spaces. We analyze relations between characteristic properties of convex bodies,
 demonstrate how $\eps$-convex bodies connect with some classical results of Convex Geometry, as Helly theorem, and find applications to geometric tomography. We introduce the notion of a circular projection and investigate the problem of determination of $\eps$-convex bodies by their projection-type images. The~results generalize corresponding stability theorems by H.\,Groemer.


\end{abstract}

\section*{Introduction}
\label{intro}

As usual, a convex body in $n$-dimensional Euclidean space $\bbr^n$ is a compact convex set with non-empty interior. Convex sets can be characterized as intersections of supporting half-spaces, see \cite{le80}. Instead of them we use bodies of some \textit{different shapes} (such as complements to balls and cylinders of "large" radius). In the present study we introduce a more general class of bodies in metric and normed linear spaces, called \textit{$\eps$-convex bodies}, in order to extend our previous results \cite{vgol}, \cite{vgol99}, \cite{Rov_2006} following to the "soft-hard" ranking of geometrical categories described in \cite{gro}.
Similar considerations were done by Reshetnyak \cite{resh56} for $\delta$-touched surfaces.
Investigation of $\eps$-convex bodies allows us to analyze relations between characteristic properties of convex bodies, to extend some classical results for bodies "close to convex", and to find applications to geometric tomography due to \cite{rgar06}, \cite{vgol}, \cite{gr87}, \cite{gr97}, etc.

It is well-known that a three-dimensional convex body is, up to translations, uniquely determined by the translates of its orthogonal projections onto all planes.
Simple examples show that this is no longer true if only "lateral projections" are permitted, that is orthogonal projections onto all planes containing a given line.
A large class of convex bodies in $\bbr^n\ (n>2)$ that are essentially determined by translates (or homothetic images) of their lateral projections is studied in \cite{gr87}, \cite{gr97}, and corresponding stability results are obtained. The orthogonal projections of different $\eps$-convex bodies onto all hyperplanes may coincide. We introduce the notion of a circular projection and investigate the problem of determination of $\eps$-convex bodies by their projection-type images. The results generalize corresponding theorems in~\cite{vgol}, \cite{gr97} and show that the class of convex bodies traditionally used in geometric tomography can be widely extended by $\eps$-convex ones with restrictions on their size (diameter,~etc).

In contrast with classical X-ray and Geometric Tomography problems \cite{rgar06}, where all measurements are connected with propagations of the signals along the straight lines, for example, in plasma tomography 
it happens so that one should consider the results of "circular" projections of the multidimensional objects under consideration.
Furthermore, in some cases one should use much more complicated than lines or circles trajectories of projections, most of the problems of photo-elasticity and seismic tomography are closely related to integrating of unknown functions, vector and tensor fields along geodesics of corresponding Riemannian metric, see for example~\cite{Shar_94}.
 Hence, the studies of $\eps$-convex bodies and their "circular" projections could be used not only in the pure theoretical domains of mathematics.

In \textbf{Section~\ref{sec:2}} we introduce basic axioms of $\eps$-convexity, study relations between new classes of bodies more complicated than convex ones (Theorems~\ref{T-Ki-eps},~\ref{T-Ki-eps4}, \ref{T-Ki-epsk}, \ref{T-Ki-epsk4}), introduce and study an $\eps$-convex hull construction, and prove a Helly type theorem (Theorem~\ref{P-epsHelly}).
In \textbf{Section~\ref{sec:3}}, devoted to applications of these classes of bodies to geometric tomography,
we study the problem: if~all circular projections (of a family similar to "lateral" ones) of two bodies in the Euclidean space are translative (homothetic) close, then these bodies are close to each other within a translative (homothetic) Hausdorff distance (Theorems~\ref{T-0eps}, \ref{T-01}).
\\


\section{$\eps$-Convex Bodies}
\label{sec:2}

In this section we introduce and study some classes of bodies in a complete metric space (in particular, the $n$-dimensional Euclidean space in the real, complex or quaternion cases, and space forms of non-zero curvature) that extend the class of convex bodies. Simply saying, the role of {separating hyperplanes} or {supporting half-spaces} for a body will play spheres or the complements to balls, resp., of radius $1/\eps$ for a given $\eps>0$. We investigate relations between these new classes of bodies, discuss an $\eps$-convex hull construction and prove corresponding analogue of Helly' theorem.

\subsection{Preliminary notions}

We remind standard notations and definitions. Let $(M,\rho)$ be a complete metric space, for example,
$n$-dimensional Euclidean space $\bbr^n$ in the real, complex or quaternion cases, $O$ its origin and $\rho(x,y)=\|x-y\|$. We~call

 $B(C,r)=\{x: \rho(x,C)\le r\}$ a closed {ball in $M$ of radius $r$ centered at~$C$},

 $\overset{o}B(C,r)=\rho(x,C)<r\}$ an open {ball in $M$ of radius $r$ centered at~$C$},

 $S(C,r)=\{x: \rho(x,C)=r\}$ a {sphere in $M$ of radius $r$ centered at~$C$}.
\newline
A point $x$ of a set $K\subset M$ is called {interior point}, if there exists $s>0$  such that  $B(x,s)$ belongs to $K$. Interior of $K$ is denoted by int\,$K$ or $\overset{o}K$.
We call $\partial{K}=K\setminus{\rm int}\,K$ the boundary of~$K$.
A set $K\subset\bbr^n$ is {closed} if int\,$(M\setminus K)=M\setminus K$.
A~diameter of a compact set $K$ is defined by $d_K=\max\{\rho(x,y): x,y\in K\}$.
Obviously, for any $x\in K$ we have $K\subseteq B(x,d_K)$.
 The distance between two non-empty sets $K,L\subset M$ is
 dist$(K,L)=\inf\limits_{x\in K, y\in L}\rho(x,y)$.
 The~Hausdorff distance between compact sets $K,L\subset M$ is defined~by
\begin{equation}\label{E-hausd}
 \delta(K,L)=\max\{\max\limits_{x\in K}{\rm dist}(x,L),\,\max\limits_{y\in L}{\rm dist}(y,K)\}.
\end{equation}
If $K\subset M$ is a non-empty set and $\eps>0$, then $K_\eps=\{x: {\rm dist}(x,K)\le\eps\}$ is called an outer parallel set to~$K$. If $K\subset\bbr^n$ then one also may use the formula $K_\eps=K +\eps B(O,1)$.
 Another definition of the Hausdorff distance is, see \cite{rgar06},
 \begin{equation}\label{E-hausd2}
 \delta(K,L)=\min\{\eps>0:\ K\subseteq L_\eps\ \ {\rm and}\ \ L\subseteq K_\eps\}.
 \end{equation}
 A subset $K$ in a metric space $(M,\rho)$
 is said to be convex if any two points of it are joined by a shortest curve in $M$ and any such shortest curve lies in $K$, see~\cite{resh4}. A convex hull $conv K$ is the smallest convex set containing $K$.
  A~convex body in $\bbr^n$ is homeomorphic to a ball, hence it is contractible and simply connected. If $K\subset\bbr^n$ is a compact body, we denote by

 $h_K(\omega)=\max\{x\cdot\omega,\ x\in K\}$ the {support function},

 $w_K(\omega)=h_K(\omega)+h_K(-\omega)>0$ the {width},

 $H_K(\omega)=\{x\in\bbr^n: x\cdot\omega=h_K(\omega)\}$ the {support hyperplane},

 $H_K^+(\omega)=\{x\in\bbr^n: x\cdot\omega\le h_K(\omega)\}$ the {support half-space}
 (containing $K$),

 $S_K(\omega)=K\cap H_K(\omega)$ the {support set}
\newline
  of $K$ in the direction~$\omega$.
 If~$S_K(\omega)$ consists of a point, it is called the {support point of $K$ in the direction $\omega$}, and $\omega$ is called a {regular direction of $K$}.

Remark that the Hausdorff distance between convex bodies $K$ and $L$ in $\bbr^n$ may be defined using the support function as follows, see \cite{rgar06}, \cite{gr87}, \cite{gr97},
\begin{equation*}
 \delta(K,L)=\max\{|h_K(\omega)-h_L(\omega)|:\ |\omega|=1\}.
\end{equation*}

\subsection{Definition and basic properties of $\eps$-convex bodies}
\label{sec:1-1}

 Recall that the hyperplane $H$ in a normed linear space $\bbr$ \textit{supports} a convex set $K$ at a point $x\in K$ if $x\in H$ and $K$ is contained in one of the half-spaces determined by $H$, \cite{nic}.
 The generalized tangency tells us that there is at least one supporting hyperplane through each boundary point of a convex body.

 The~basic separation theorem tells us that
 {if a point $x$ is disjoint from a (compact) convex body $K$ in a normed linear space $E$ then there exists a closed hyperplane that separates strictly $K$ and~$x$}.

Starting from well-known characteristic properties of {convex bodies}
in terms of supporting and separating hyperplanes or half-spaces, we introduce some classes of compact bodies in a metric space $(M,\rho)$ generalizing convex ones.

A~ball $B$ will be called an \textit{outer support ball} of a body $K$ if it doesn't intersect int\,$K$ and intersects $K$ at its boundary points. Denote by $\mathcal{C}(K,\eps)$ a set of centers of all outer support balls $B(C,1/\eps)$ of a body $K$ in a metric space $(M,\rho)$.

\begin{defi}\label{D-eps-conv}\rm
A body $K\subset M$ is called $\eps$-\textit{convex} of a class $\mathcal{K}^\eps_i$
(for some $\eps>0$)~if

\hskip-1mm
$\mathcal{K}^\eps_1$:
{any point $x\in\partial K$ belongs to an outer support ball (of $K$)
of radius~$1/\eps$},

\vskip.5mm
\hskip-1mm
$\mathcal{K}^\eps_2$:
{any point $x\notin K$ belongs to a ball $B$ of radius $1/\eps$ such that
int\/$K\cap B{=}\emptyset$},

\vskip.5mm
\hskip-1mm
$\mathcal{K}^\eps_3$: {any outer support ball $B(x,r)$ of $K$ with $r<1/\eps$ belongs to
an outer}

{support ball (of $K$) of radius $1/\eps$}.

\noindent
A~connected boundary component of an $\eps$-convex body in $\bbr^n$ will be called an $\eps$-\textit{convex hypersurface} ($\eps$-{convex curve} if $n=2$) of a certain class listed above.
\end{defi}

\begin{prop}\label{P-Ki-eps4}
The class $\mathcal{K}^{\eps}_3$ of bodies in complete metric spaces coincides with the class defined by the following weaker condition:

\vskip.5mm
\hskip-1mm
$\mathcal{K}^\eps_4$: {any ball $B(y,r)$ of radius $r<1/\eps$ disjoint from $K$ belongs to a ball of}

 {radius $1/\eps$ that does not intersect $int\,K$}.
\end{prop}

\textbf{Proof}. We need to show the inclusions.

$\mathcal{K}^\eps_3\subseteq\mathcal{K}^\eps_4$.
Given $K\in\mathcal{K}^\eps_3$ and a ball $B(y,r)$ of radius $r<1/\eps$ disjoint from $K$,
consider an outer support ball $B(y,\rho)$ of $K$, where $\rho={\rm dist}(y,K)$.
If~$\rho\ge1/\eps$ then $B(y,1/\eps)$ is a desired ball.
 If $\rho<1/\eps$ then $B(y,r)$ is contained in an outer support ball $B(y,\rho)$.
 Due to $(\mathcal{K}^\eps_3)$, the ball $B(y,\rho)$ is contained in some outer support ball $B(z,1/\eps)$ of $K$.

$\mathcal{K}^\eps_4\subseteq\mathcal{K}^\eps_3$.
 Given $K\in\mathcal{K}^\eps_4$ and a ball $B(y,r)$ of radius $r<1/\eps$ supporting $K$, consider a sequence of balls $B(y,r-\frac1m)$ disjoint from $K$. Due to ($\mathcal{K}^\eps_4$), each ball $B(y,r-\frac1m)$ belongs to a ball
 $B(y_m,1/\eps)$ disjoint from int$K$.
 There is a subsequence $m_s\ (s\to\infty)$ such that $y_{m_s}{\to}\, y_0$.
 Then the ball $B_0=B(y_0,1/\eps)$ is disjoint from int$K$ and contains $B(y,r)$. Moreover, $B_0$ is supporting to~$K$.$\,\square$

\begin{rem}\rm (a) One may verify (applying just the set theory arguments and Proposition~\ref{P-Ki-eps4}) that if a body $K$ is the intersection of (connected) bodies of a class $\mathcal{K}^\eps_2$ (or $\mathcal{K}^{\eps}_1$)
then $K$ also belongs to $\mathcal{K}^\eps_2$ (resp., $\mathcal{K}^{\eps}_1$).
 This intersection can be disconnected, see examples in what follows. Hence bodies in $\bbr^n$ of a class $\mathcal{K}^\eps_2$ can be represented as the sets of solutions $x$ to a systems of weak quadratic inequalities $(x-C_\alpha)^2\ge1/\eps^2,\ \alpha\in I$.

(b) If $M=\bbr^n$, one may replace the balls of radius $1/\eps$ in Definition~\ref{D-eps-conv} by half-spaces
 and obtain (when the bodies are assumed connected) the class $\mathcal{K}^0$ of convex bodies in $\bbr^n$. Obviously, for $n=1$ the classes $\mathcal{K}^\eps_i$ coincide.

(c) Similar classes of $\eps$-convex bodies can be defined in complex Euclidean spaces $\bbc^n$, where the complex hyperplanes have real codimension 2 and do not divide the space. Boundaries of balls and cylinders over complex (quaternion) affine subspaces are hypersurfaces of real codimension~1, hence they can be used in definition of classes $\eps$-convex bodies in $\bbc^n\ (\bbh^n)$ as~well.
\end{rem}

\begin{exam}\rm
(a) Simple examples of $\eps$-convex bodies can be obtained using two balls in real Euclidean spaces $\bbr^n$, one of them has radius at least $1/\eps$. A~homeomorphic to a ball body $K=B(O,r)\setminus\overset{o}B(C,1/\eps)$ (for any point $C$ with $|OC-1/\eps|<r$) belongs to a class $\mathcal{K}^\eps_i$.
A~body $K=B(O,r)\setminus\overset{o}B(O,1/\eps)=\{x\in\bbr^n:\ 1/\eps\le\|x\|\le r\}$
with $r>1/\eps$ belongs to a class $\mathcal{K}^\eps_i$, its boundary has two components.
Slight modification leads to a non-concentric version of this ring, say:
$K=B(O,N/\eps)\setminus\overset{o}B(O_1,1/\eps)$, where $|OO_1|<(N-1)/\eps$.
 Remark that $I^n\setminus B(C,1/\eps)$, where is $I^n$ is a "huge" cube (with side $2/\eps$) and center at $O$, belongs to a class $\mathcal{K}^\eps_i$ as well.

(b) A non-convex polyhedra in $\bbr^n$ does not belong to $\mathcal{K}^\eps_i$ for any $\eps>0$ and $i=1,2,3$. Simple examples are the quadrangle $OABC\subset\bbr^2$ where $A(1,0),B(0,1)$, $C(1/4,1/4)$, and cylinders $OABC\times[-1,1]^{n-2}\subset\bbr^n$.

(c) Let a $C^2$-regular hypersurface $\partial K\subset\bbr^n$ bounds a body $K\in\mathcal{K}^\eps_1$, and let $\bn$ be a unit normal to $\partial K$ directed inside. Denote by $\kappa(\bn,\ba)$ the normal curvature of $\partial K$ with respect to $\bn$ in the direction $\ba$. Set $\kappa_{\partial K}=\min_{\ba}\kappa(\bn,\ba)$.
Then $\kappa_{\partial K}\ge-\eps$. In particular, the curvature of an $\eps$-convex curve satisfies the inequality $\kappa\ge-\eps$. For example, astroid $\gamma: x^{2/3}+y^{2/3}=1$
is a curve  of a class $\mathcal{K}^{2/3}_1$. Note that $d_\gamma=2>1/(2\eps)$, i.e., astroid bounds a "large" body in $\bbr^2$.

(d) An example of a "small" (see also Section~\ref{sec:1.3}) disconnected body in $\mathcal{K}^\eps_i\ (i=1,2)$ can be obtained using three discs in $\bbr^2$. Namely, let $K'=B(O,r)\setminus(B(C,1/\eps)\cup B(-C,1/\eps))\subset\bbr^2$, where $C=(r/2,0)$ and $0<r<0.8/\eps$. Then a cylindrical body $K=K'\times[0,r/2]^{n-2}\subset\bbr^n$ has two components and belongs to $\mathcal{K}^{\eps}_i$. A~"large" disconnected body in $\mathcal{K}^\eps_i\ (i=3,4)$ is represented by a union of two balls $K=B(a,r)\cup B(-a,r)$ in $\bbr^n$, where $0<r<\|a\|$ and $\|a\|-r>1/\eps$.

(e) The bodies
$
 K_1=\{y\ge0.1\cos x,\,|x|\le2\pi,\,y\le1\},\
 K_2=\{y\le-0.1\cos x,$ $|x|\le2\pi,\,y\ge-1\}
$
are homeomorphic to a ball in $\bbr^2$ and belong to a class~$\mathcal{K}^{0.1}_2$.
Their intersection $K=K_1\cap K_2$ has two convex components given by $\{0.1\cos x\le y\le-0.1\cos x,\,|x|\le2\pi\}$. Two bodies (translates) $K_1+(0,0.1)$ and $K_2-(0,0.1)$ belong to a class $\mathcal{K}^{0.1}_2$ and they intersect at two distinct points $(\pm\pi,0)$.
The neighborhoods of a segment $[A,B]$ and any arc $AB$ with small curvature are $\eps$-convex and their intersection consist of neighborhoods of two points $\{A,B\}$.

(f) We will build a simply connected not contractible body $K\subset\bbr^n\ (n\ge5)$ of a class $\mathcal{K}^{\eps}_2$. Denote by $S^1(O,r)\subset\bbr^n\ (n\ge5)$ a circle in the plane $P_1=\{x_3=\ldots=x_n=0\}$ of radius $r=1/\eps-\eps$ for small $\eps$. Let $L$ be a union of all balls $B(M, 1/\eps)$, where $M\in S^1(O,r)$. Obviously, a body $K=B(O,3)\setminus L$ belongs to a class $\mathcal{K}^{\eps}_2$. We will show that $K$ is homotopy equivalent to a standard sphere $S^{n-3}$, from this (since $n-3\ge2$) will follow the desired property.
Set $K_1=B_1\setminus L$, where $B_1=B(O,3)\cap P_1^\perp$ a 3-dimensional ball.
From the analytic representation $K_1=\{x\in\bbr^n:\ x_1=x_2=0,\ \sqrt{2-\eps^2}\le\|x\|\le3\}$
it follows that $K_1$ is homotopy equivalent to a sphere $S^{n-3}$.
Moreover, one can show that $K_1$ is a retract of $K$ (the retraction is organized along segments orthogonal to $P_1$). Hence $K$ is a desired body. Note that for $n=4$ this construction leads to a non simply connected body $K$ (namely, homotopy equivalent to a circle) of a class $\mathcal{K}^{\eps}_2$.
Further examples of $\eps$-convex bodies are presented in the sequel.
\end{exam}

\begin{theo}\label{T-Ki-eps} The $\eps$-convex bodies in complete metric spaces satisfy the following strong inclusions: $\mathcal{K}^\eps_3\subset\mathcal{K}^\eps_2\subset\mathcal{K}^\eps_1$.
\end{theo}

\textbf{Proof}. (a) $\mathcal{K}^\eps_2\subset\mathcal{K}^\eps_1$.
Given $K\in\mathcal{K}^\eps_2$ and $x\in\partial K$, consider a sequence $x_n\not\in K$ such that $x_n\to x$. Due to $(\mathcal{K}^\eps_2)$, for each $n$ there is a ball $B(C_n,1/\eps)$ containing $x_n$ such that int\,$K\cap B(C_n,1/\eps)=\emptyset$.
The sequences $C_n$ is bounded, hence we may assume that $C_n\to C$. From this it follows that $1/\eps\ge\lim\nolimits_{n\to\infty}\rho(x_n,C_n)=\rho(x,C)$.
Then $x\in B(C,1/\eps)$ and $B(C,1/\eps)\cap{\rm int}\,K=\emptyset$ as required.

\begin{figure*}[ht]
\begin{center}
\begin{tabular*}{\textwidth}{@{\hskip3mm}l@{\extracolsep{\fill}}r@{\hskip3mm}}
 \subfigure[\label{F-cond-ia} $K'\in\mathcal{K}^\eps_1\setminus\mathcal{K}^\eps_2$.]{%
\includegraphics[scale=  0.23,angle=0  ]{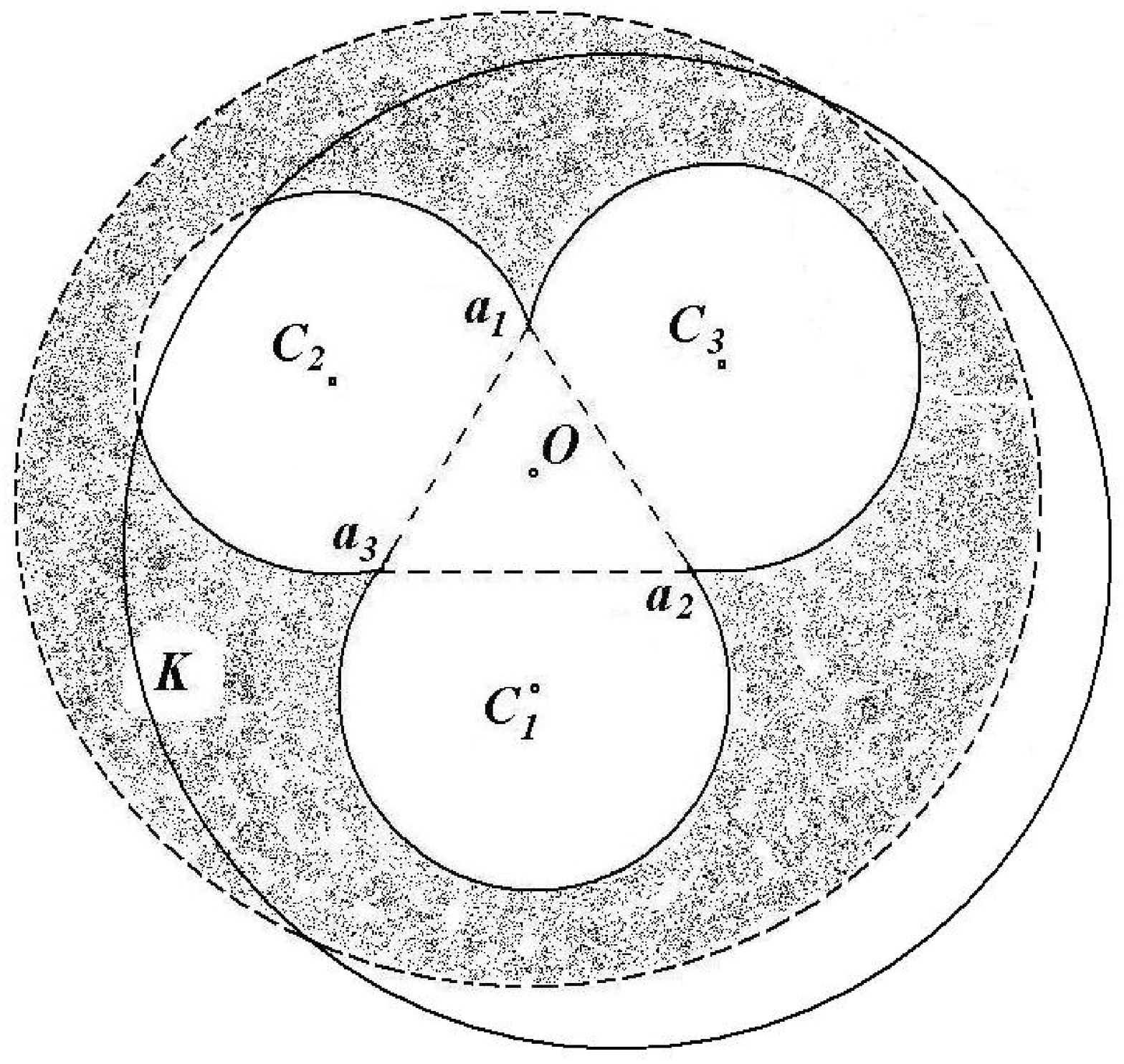}} &
 \subfigure[\label{F-exam2b} $K\in\mathcal{K}^\eps_2\setminus\mathcal{K}^\eps_3$.]{%
\includegraphics[scale=  0.45,angle=0  ]{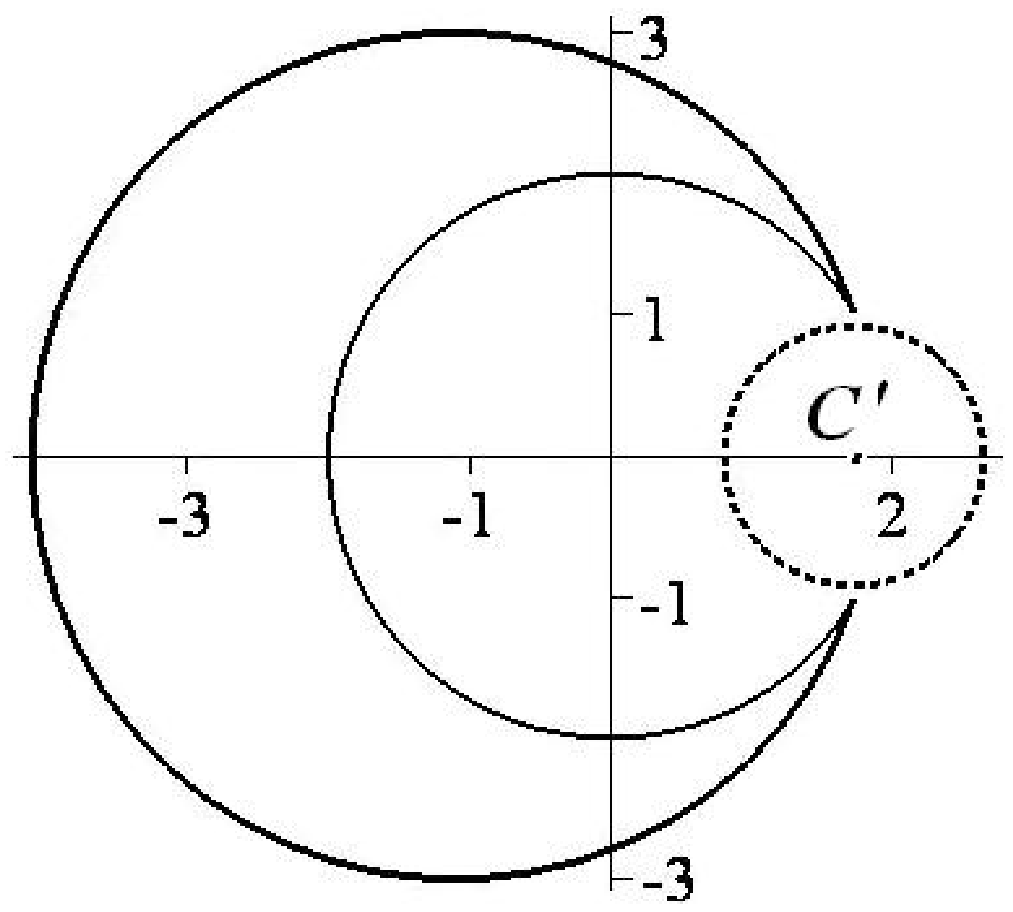}}
\end{tabular*}
\caption{\small }
 \label{F-08-011}
\end{center}
\end{figure*}

Let us show that in $\bbr^n\ (n\ge2)$ there is a homeomorphic to a ball body $K\in\mathcal{K}^\eps_1\setminus\mathcal{K}^\eps_2$. Cut from a ball $B(O,4r)$ a regular simplex $\Delta_n=conv\{a_1,\dots a_{n+1}\}$, i.e., a~convex hull of $n+1$ equally distanced points (an equilateral triangle when $n=2$) with center at $O$, and $n+1$ balls $B(C_1,r),\dots,B(C_{n+1},r)$ such that
$\partial B(C_i,r)$ contains $\{a_1,\dots a_{n+1}\}\setminus\{a_i\}$. One may select
${\rm dist}(a_i,a_j)<2\,r$ so that $O$ is not contained in any of these balls, see Fig.~\ref{F-cond-ia} for $n=2$. The body
\begin{equation}\label{E-K23}
\begin{array}{c}
 K=B(O,4r)\setminus(\Delta_n\cup\bigcup\nolimits_i B(C_i,r))\subset\bbr^n
\end{array}
\end{equation}
satisfies $(\mathcal{K}^\eps_1)$ for $r=1/\eps$, but $(\mathcal{K}^\eps_2)$ is not satisfied for a point $O$. If we slightly move  a ball $B(O,4r)$ in the direction $C_2O$, then we will obtain a homeomorphic to a~ball body $K'$ with similar properties, see Fig.~\ref{F-cond-ia}.

\vskip.5mm
 (b) $\mathcal{K}^\eps_3\subset\mathcal{K}^\eps_2$.
 Let $K\in\mathcal{K}^\eps_3$. Suppose that $x\not\in K$. If ${\rm dist}(x,K)\ge1/\eps$ then $B(x,1/\eps)\cap{\rm int}\,K=\emptyset$, and $B(x,1/\eps)$ is a desired ball.
 If dist$(x,K)<1/\eps$ then dist$(x,K)\le 1/\eps-r$ for any $r\in(0,1/\eps-{\rm dist}(x,K))$.
 By condition $(\mathcal{K}^\eps_3)$, there is an outer support ball $B(C,1/\eps)\supset B(x,r)\ni x$, hence $K\in\mathcal{K}^\eps_2$.

Let us show that in $\bbr^n\ (n\ge3)$ there is a homeomorphic to $n$-ball body $K\in\mathcal{K}^\eps_2\setminus\mathcal{K}^\eps_3$. Suppose $n=3$. Consider two discs $K_1=\{(x+1.1)^2+y^2\le9,\ z=0\}$ and $K_2=\{x^2+y^2\le4,\ z=0\}$, Fig.~\ref{F-exam2b}.
The shape of a cylindrical set $W=(K_1\setminus K_2){\times}[-0.1,0.1]\subset\bbr^3$ is as of a letter "c". Two balls, $B_1=B(C,1/\eps)$ and $B_2 =B(-C,1/\eps)$, where $C=(0,0,\sqrt{1/\eps^2{-}4})$ and $\eps>0$ is small, have a common circle $\gamma=\{x^2{+}y^2=4,\ z=0\}$. A~"small" body $K=W\setminus(B_1\cup B_2)$ is homeomorphic to a ball, belongs to~$\mathcal{K}^\eps_2$, and $B_1,B_2$ are outer support balls of radius $1/\eps$ of $K$. The balls $B(O,r)$, where $0<r<2$, and $B(C',r')$, where $C'(x',0,0),\ x'=(5-1.1^2)/2.2\approx1.72$, $0<r'\le0.5$, do not intersect~$K$. Obviously, no outer support ball of $K$ of radius $1/\eps$ contains $B(O,r)$ or $B(C',r')$. Hence a body $K$ doesn't belong to $\mathcal{K}^\eps_3$. Remark that condition $(\mathcal{K}^\eps_3)$ is not valid for a point $C'$.
 One may modify this example using $K_1=\{x^2+y^2\le9,\ z=0\}$ to obtain a non simply connected $\eps$-convex body $K\in\mathcal{K}^\eps_2\setminus\mathcal{K}^\eps_3$, shaped as a letter "o". Remark that similar bodies $K$ and $K'$ exist in $\bbr^n$ for any $n\ge3$.$\,\square$

\vskip1.5mm
The next theorem extends the Motzkin's characterization of convex bodies in $\bbr^n$ \cite{mot35} (see also Theorem~9.3 in~\cite{le80}) for a larger class of bodies.

\begin{theo}\label{T-Ki-eps4}
The class $\mathcal{K}^{\eps}_3$ of bodies in $\bbr^n$ coincides with the class defined by the following condition:

\vskip.5mm
\hskip-3mm
$\mathcal{K}^\eps_{5}$: {for any point $z\in M\setminus K$ such that $\rho={\rm dist}(z, K)<1/\eps$ the intersection}

{$B(z,\rho)\cap K$ consists of exactly one point}.
\end{theo}

\begin{figure}[ht]
\begin{center}
\includegraphics[scale=.3,angle= 0,clip=true,draft=false]{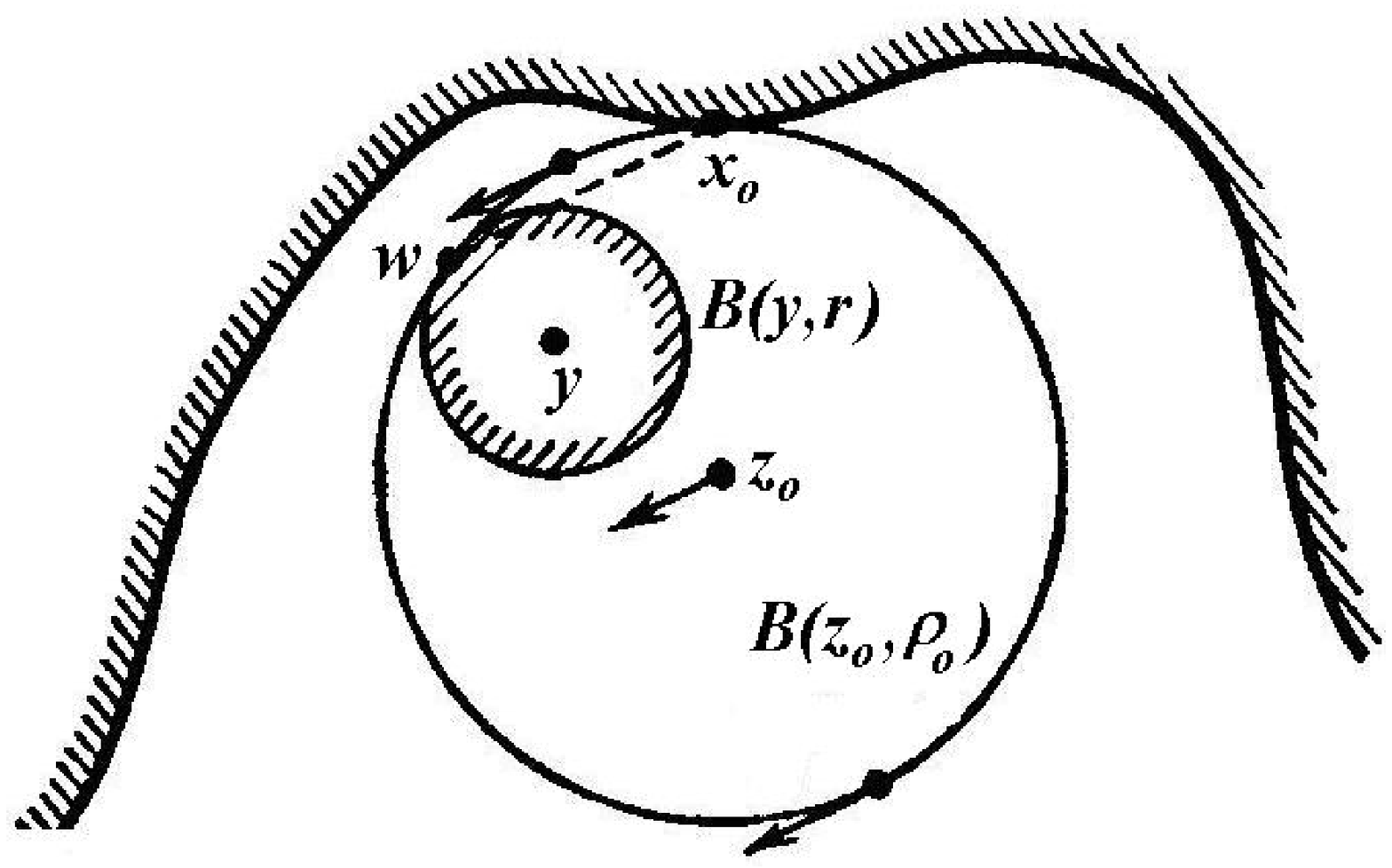}
\caption{\small Proof: $\mathcal{K}^\eps_{5}\subseteq\mathcal{K}^\eps_3$.}
 \label{F-Lei28}
\end{center}
\end{figure}

\textbf{Proof}. $\mathcal{K}^\eps_{5}\subseteq\mathcal{K}^\eps_3$.
Let $K\in\mathcal{K}^\eps_{5}$. Suppose an opposite, that there is an outer support ball $B(y,r)$ of $K$ of radius $r<1/\eps$ that is not contained in any outer support ball of radius $1/\eps$.
Denote by $\Omega$ a set of all outer support balls $B(z,\rho)\subset\bbr^n$ of $K$ such that
\begin{equation}\label{E-L142}
 \rho\le1/\eps,\qquad
 B(y,r)\subset B(z,\rho).
\end{equation}
 Obviously, $B(y,{\rm dist}(y,K))\in\Omega$, hence $\Omega\not=\emptyset$.
 From (\ref{E-L142}) it follows (since $K$ compact and the radii $\rho\le1/\eps$) that the centers $z$ of these balls form a bounded set.
 Hence there are $\rho_0=\sup\{\rho:\ B(z,\rho)\in\Omega\}\le1/\eps$
 and a sequence $\{B(z_i,\rho_i)\}_{i\in\bbn}$ such that $\lim\limits_{i\to\infty}\rho_i=\rho_0$.
 We can assume that $\lim\limits_{i\to\infty} z_i=z_0$. For $B(z_0,\rho_0)$ holds
  \begin{equation}\label{E-L144}
  B(y,r)\subset B(z_0,\rho_0),\quad
  B(z_0,\rho_0)\cap\,int\,K=\emptyset,\quad
  B(z_0,\rho_0)\cap K\not=\emptyset,
 \end{equation}
 thus $B(z_0,\rho_0)$ is the maximal ball from $\Omega$.
 By our assumption, $\rho_0<1/\eps$, otherwise $B(z_0,1/\eps)$ is an outer support ball containing $B(y,r)$.
 Denote by $x_0$ a (unique due $(\mathcal{K}^\eps_5)$) nearest to $z_0$ point of $K$, i.e., $\rho_0=\rho(x_0,z_0)$.
 Since $B(y,r)\cap K=\emptyset$, we have $r<\rho_0$. Hence the boundary sphere $S(z_0,\rho_0)$ of a ball $B(z_0,\rho_0)$ has at most one common point with $B(y,r)$.
  If there are no common points, we move $S(z_0,\rho_0)$ onto small distance in the direction $z_0-x_0$, and its image $S(z'_0,\rho_0)$ will have the property
 \begin{equation}\label{E-Lprof}
  B(y,r)\subset B(z'_0,\rho_0),\quad
  B(z'_0,\rho_0)\cap K=\emptyset.
 \end{equation}
 Then $B(z'_0,\rho')\in\Omega$ for some $\rho'={\rm dist}(z'_0,K)>\rho_0$
 that contradicts to maximality of $B(z_0,\rho_0)$, see Fig.~\ref{F-Lei28}.
 If there is one common point $w$, we move $S(z_0,\rho_0)$ onto small distance in the direction $w-x_0$,
 and its image $S(z'_0,\rho_0)$ will have the property (\ref{E-Lprof}), again a contradiction.

\vskip.5mm
$\mathcal{K}^\eps_3\subseteq\mathcal{K}^\eps_{5}$.
 Let $K\in\mathcal{K}^\eps_3$. Suppose an opposite, that there is a point $z\in\bbr^n\setminus K$ such that $\rho={\rm dist}(z,K)<1/\eps$ and there are two different points $x_1,x_2\in K$ with the property $\rho(x_1,z)=\rho(x_2,z)=\rho$. Then $B(z,\rho)$ is an outer support ball of $K$.
 Due to $(\mathcal{K}^\eps_3)$, $B(z,\rho)$ is contained in an outer support ball $B(C,1/\eps)$ of $K$. Two spheres $S(z,\rho)$ and $S(C,1/\eps)$ have at most one common point. Since $x_1,x_2\in S(z,\rho)$, one may assume that $x_1\not\in S(C,1/\eps)$.
 Hence $x_1$ is an inner point of $B(C,1/\eps)$ that is a~contradiction.$\,\square$

\begin{rem}\rm The claim of Theorem~\ref{T-Ki-eps4} is wrong for the linear space $(\mathbb{R}^2, \|\ \|_\infty)$ with the norm $\|(x,y)\|_\infty{=}\max\{|x|,|y|\}$. In this case, a ball $B(0,1)$ is presented by a unit square, that is not strictly convex. Hence the convex bodies in normed linear spaces may not satisfy condition ($\mathcal{K}^\eps_{5}$). Nevertheless, Theorem~\ref{T-Ki-eps4} may be extended for all normed linear spaces with strictly convex unit balls.
\end{rem}

\subsection{Helly type theorem for "small" $\eps$-convex bodies}
\label{sec:1.3}

An $\eps$-convex body $K$ will be called "small" if $\eps d_K\le1$, and $K$ will be called "large" if $\eps d_K>1$. Next proposition illustrates this showing that "small" $\eps$-convex bodies don't contain holes inside of their interior.

\begin{prop}\label{P-0}
If a "small" body $K\in\mathcal{K}^\eps_1$ is connected, then its boundary $\partial K$ is connected.
\end{prop}

\textbf{Proof}. Let $K$ be connected and its boundary $\partial K$ is not, then one of components of $\partial K$ i.e., $\partial_1 K$, separates $K$ from "infinity"
(or from $\bbr^n\setminus conv K_\eps$).
Now, $\partial_1 K$ is a boundary of some body $L$ such that $d_L=d_{\partial K}=d_K$.
Let $\partial_2 K\neq\partial_1 K$ be another component of $\partial K$ and let $x_1\in\partial_2 K$.
If $K\in\mathcal{K}^\eps_1$ then $x_1$ should be contained in $L$ and in some outer support ball $B(C,1/\eps)$. So, $B(C,1/\eps)\subset L$ and $d_L\ge1/\eps$, a contradiction.
Remark that for $K\in\mathcal{K}^\eps_i\ (i=2,3,4)$ the claim is also true due to Theorem~\ref{T-Ki-eps}.$\,\square$

\vskip1.5mm
From the definition it follows that if $L\subseteq K$ be compact sets in $\bbr^n$
then $h_L(\omega)\le h_K(\omega)$ for all $\omega$, or simply $h_L\le h_K$.
The opposite holds when $K$ is a convex set, see \cite{rgar06}, p.\,16.

\begin{prop}\label{P-ed2}
Let a body $K\subset\bbr^n$ belongs to a class $\mathcal{K}^\eps_3$, and $d_K\le\frac1{2\eps}$, and let $L\subset\bbr^n$ be a compact set. If $\,h_L\le h_K$ then $L\subseteq K_{\eps'}$, where $\eps'=\eps d^2_K/2$.
\end{prop}

\begin{figure}[ht]
\begin{center}
\includegraphics[scale=.45,angle= 0,clip=true,draft=false]{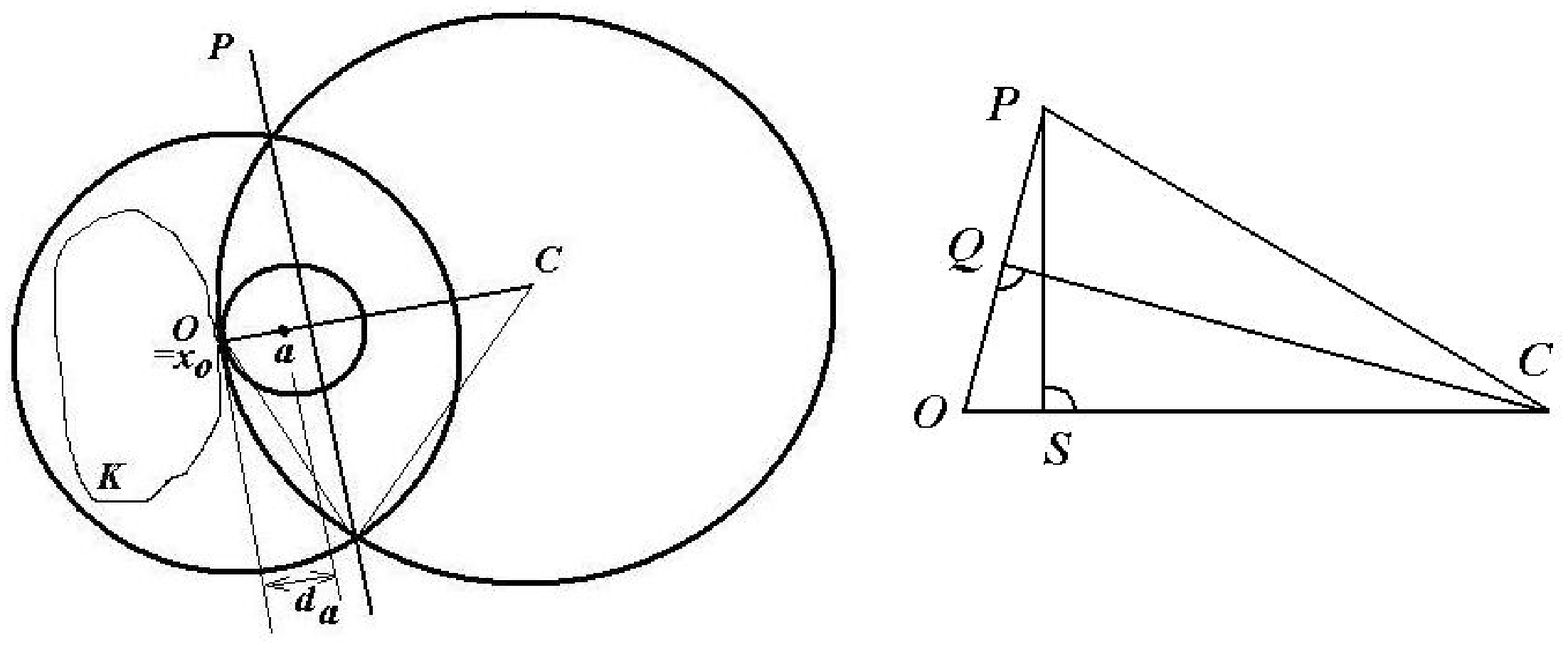}
\caption{\small Support function of $K\in\mathcal{K}^\eps_3$.}
 \label{F-prop21}
\end{center}
\end{figure}

\textbf{Proof}. Let there is a point $a\in L$ such that $d_a={\rm dist}(a,K)>0$ (otherwise, $L\subseteq K$ and there is nothing to prove). Denote by $x_0$ the (unique) nearest to $a$ point of $K$.
 Since $h_{K-x_0}(\omega)=h_{K}(\omega)-x_0\cdot\omega\le h_{L}(\omega)-x_0\cdot\omega=h_{L-x_0}(\omega)$,
 we may assume $O=x_0$, hence $K\subset B(O,d_K)$ and $d_a=\|a\|$, see Fig.~\ref{F-prop21}.

\noindent
If $\|a\|\ge d_K$ then $h_L(a/\|a\|)\ge\|a\|\ge d_K>h_K(a/\|a\|)$ -- a contradiction.
Hence, $\|a\|<d_K\le1/(2\,\eps)$. Set $d'_a=d_a-\alpha$, where $\alpha>0$ is sufficiently small.
 By Theorem~\ref{T-Ki-eps}, there is a (outer support to $K$) ball $B(C',1/\eps)$ separating $K$ from a ball $B(a,d'_a)$. Obviously, $C'\to C$ for $\alpha\to0$. Hence $O\in B(C,1/\eps)\cap K$.
 Denote by $P(a)$ a hyperplane (orthogonal to $a$) through the intersection of the spheres $S(C,1/\eps)$ and $S(O,d_K)$. From elementary geometry we obtain
 $OS/OP=PQ/PC$, where $OP=2\,PQ=d_K,\,PC=1/\eps$ and $OS>\|a\|$, Fig.~\ref{F-prop21},
 hence ${\rm dist}(O,P(\omega))=\eps d^2_K/2$.
 Remark that $h_K(\omega)\le{\rm dist}(O,P(\omega))$ and $d_a\le h_L(\omega)$.
From this and $h_L(\omega)\le h_K(\omega)$ it follows that $d_a\le\eps d^2_K/2$.
$\,\square$

\begin{rem}\rm The condition $K\in\mathcal{K}^\eps_3$ in Proposition~\ref{P-ed2}
can not be replaced by $K\in\mathcal{K}^\eps_2$ for $n\ge3$.
 Really, let $K\in\mathcal{K}^\eps_2\setminus\mathcal{K}^\eps_3$ be a homeomorphic to a ball
body in $\bbr^3$, presented in the proof of Theorem~\ref{T-Ki-eps}.
We have $d_K\le3$. Set $L=K\cup\{C'\}$, where $C'\approx(1.72,0,0)$. Obviously, the equality $h_L=h_K$ holds.
 The distance from any of two singular points on $K$ from $x$-axis is $h\approx1.01$.
 Hence $\delta(L,K)\ge{\rm dist}(C',K)\ge h>1$.
 If we take $\eps<0.2$ then $d_K<\frac1{2\eps}$ and $\eps d^2_K/2\le 0.9$.
 In this case $L$ does not belong to a class $K_{\eps'}$ for $\eps'=\eps d^2_K/2$.
\end{rem}

From Proposition~\ref{P-ed2} it follows

\begin{cor}
Suppose that bodies $L,K\subset\bbr^n$ belong to a class $\mathcal{K}^\eps_3$,
and $d=\max\{d_L, d_K\}\le1/(2\,\eps)$ holds. If $h_L=h_K$ then $\delta(L,K)\le\eps d^2/2$.
\end{cor}

Remark that $h_K=h_{conv K}$. From Proposition~\ref{P-ed2} it follows

\begin{cor}\label{C-22}
 Suppose that a body $K\subset\bbr^n$ belongs to a class $\mathcal{K}^\eps_3$, and $d_K\le1/(2\,\eps)$ holds. Then $\delta(K, conv K)\le\eps d^2_K/2$.
\end{cor}

Next theorem generalizes well-known Helly theorem ($\eps=0$).

\begin{theo}[$\eps$-Helly]\label{P-epsHelly} Let $K^1,\dots K^m\ (m>n)$ be bodies in $\bbr^n$ of diameter $d_{K^i}\le d<1/(2\eps)\ (i\le m)$ of a class $\mathcal{K}^\eps_3$.
Suppose that for any $K^{i_1},\dots K^{i_{n{+}1}}$ of this family there is $x'\in\bbr^n$ such that {\rm dist}$(x',K^{i_j})\le\eps'$ for $j\le n{+}1$ ($\eps'=0$ if their intersection is non-empty). Then {\rm dist}$(x,K^i)\le\eps'+\eps\/(m{-}n{-}1)\,d^2/2$ for some $x\in\bbr^n$
and all $i\le m$.
\end{theo}

\textbf{Proof}. We apply induction for a number $s$ of bodies.
For $s=n+1$ the claim is true due the conditions. Suppose that the claim is true for any family of $n+1\le s<m$ bodies satisfying conditions of theorem. Consider a family of $s=m\ge n+2$ bodies of a class $\mathcal{K}^\eps_3$ satisfying these conditions. By induction hypothesis, for each $i\le m$ there is $a_i\,{\in}\bigcap_{j\not=i} K^j_{\eps'}$ where $\eps''=\eps'+\eps\/(m{-}n{-}2)\,d^2/2$.
 A set $A=\{a_1,\dots a_m\}$ consists of at least $n+2$ points, and by Radon theorem (see, for example, \cite{le80}) it can be divided into two non-intersecting subsets, $A=A'\cup A''$, such that $conv A'\cap conv A''\not=\emptyset$.  We can re-order indices so that
 $A'=\{a_1,\dots a_k\}$ is contained in $K^i_{\eps'}$ for all $i>k$,
 and
 $A''=\{a_{k+1},\dots a_m\}$ is contained in $K^i_{\eps''}$ for all $i\le k$.
 Since $conv A'\cap conv A''\not=\emptyset$, there is a convex combination with two representations, $a_\lambda=\sum_{i\le k}\lambda'_i a_i=\sum_{i>k}\lambda'_i a_i$,
 where $\sum_{i\le k}\lambda'_i=\sum_{i>k}\lambda'_i=1$ and $\lambda'_i\ge0$.
   In view of dist$(a_i, K^j)\le\eps''$ for $i\le k,\,j>k$, there are $a'_{ij}\in K^j$ such that $\|a_i-a'_{ij}\|\le\eps''$ for $i\le k,\,j>k$.
 By Corollary~\ref{C-22}, dist$(a'_{j\lambda}, K^j)\le\eps''\/d^2/2$ for $j>k$,
 where $a'_{j\lambda}=\sum_{i\le k}\lambda'_i a'_{ij}$.
 Remark that
 $$
 \begin{array}{c}
 \|a_\lambda-a'_{j\lambda}\|=\|\sum_{i\le k}\lambda'_i(a_i-a'_{ij})\|
 \le\sum_{i\le k}\lambda'_i\|a_i-a'_{ij}\|\le\eps''\sum_{i\le k}\lambda'_i=\eps''.
 \end{array}
 $$
 Hence, for $j>k$ we obtain
 $$
 {\rm dist}(a_\lambda, K^j)\le\|a_\lambda-a'_{j\lambda}\|+{\rm dist}(a'_{j\lambda}, K^j)
 \le\eps''+\eps\/d^2/2=\eps'+\eps\/(m{-}n{-}1)\,d^2/2.
 $$
 Similarly, dist$(a_\lambda, K^j)\le\eps'+\eps\/(m{-}n{-}1)\,d^2/2$ for $j\le k$.
 From above follows that $x=a_\lambda$ is a desired point.$\,\square$


\subsection{$\eps$-Convex hull}

We define and examine here the $\eps$-convex hull construction.

\begin{defi}\label{D-eps-convM}\rm
Let $K$ be a compact set in a metric space $(M,\rho)$ and let $conv K$ be a body.
The intersection of all bodies $K'\subset M$ of a class $\mathcal{K}^\eps_2$ containing~$K$
is called an \textit{$\eps$-convex hull} of $K$ and is denoted by $conv_\eps K$.
\end{defi}

Obviously,  $conv_\eps K$ belongs to $conv K$.
For subsets in Euclidean spaces Definition~\ref{D-eps-convM} will be completed by details.
Denote by aff\,$K$ the \textit{affine hull} of a set $K\subset\bbr^n$, a minimal affine subspace containing~$K$. We call dim\,aff\,$K\le n$ the \textit{affine dimension} of $K$.
In other words, dim\,aff\,$K$+1 is the maximal number of  affine independent points in $K$.

\begin{defi}\label{D-eps-convh}\rm
Let $K\subset\bbr^n$ be a compact set, and dim\,aff\,$K=m$. Identify ${\rm aff}\,K=\bbr^m$ by any isomorphism. The intersection of all bodies $K'\subset\bbr^m$ of a class $\mathcal{K}^\eps_2$ containing $K$ is called an \textit{$\eps$-convex hull} of $K$ and is denoted by $conv_\eps K$.
\end{defi}

 Obviously, $conv_\eps K$ belongs to $\bbr^m$.  Since a convex body is $\eps$-convex for all $\eps>0$ we also have $conv_\eps K\subset conv K$.

\begin{rem}\label{R-convhull}\rm
 Another definition of $conv_\eps K$, when aff\,$K=\bbr^m$, is as follows:
\begin{equation*}
\begin{array}{c}
 conv_\eps K=\bbr^m\setminus\cup B^m(C,1/\eps),\quad
 {\rm where}\quad
 int\,K\cap{B^m}(C,1/\eps)=\emptyset.
\end{array}
\end{equation*}
This means that if $int{B^m}(C,1/\eps)\cap K=\emptyset$ then $int{B^m}(C,1/\eps)\cap conv_\eps K=\emptyset$.
 \end{rem}

We will briefly \textbf{prove} the equivalence of two definitions
of an $\eps$-convex hull (see Definition~\ref{D-eps-convh} and Remark \ref{R-convhull}):

 A) $conv_\eps K=\bigcap_*K_*$, where $K_*$ are all $\eps$-convex bodies containing $K$;

 B) $conv_\eps K=E^m\setminus\bigcup B^m(C,1/\eps)$, where $int\,K\bigcap{B^m}(C,1/\eps)=\emptyset$.

\vskip1mm
\noindent	
In order to prove  $(A)\Rightarrow(B)$ consider $y\in\bigcap_*K_*$
(union of all $K_*$ such that it is contained in some $B^m(C_0,1/\eps)$ such that
$int\,K\bigcap{B^m}(C_0,1/\eps)=\emptyset$). Let $L_*=K_*\bigcap(E^m\setminus B^m (C_0,1/\eps)$, (intersection of bodies of a class $\mathcal{K}^\eps_2$ belongs to $\mathcal{K}^\eps_2$) we have $K\subset L_*$ for all these $L_*$ and $y\notin L_*$ for all these $L_*$, a contradiction. In a similar way one can prove that $(B)\Rightarrow (A)$.
$\,\square$

\begin{prop}\label{P-23} Let $K\subset\bbr^n$ be a compact set, and dim\,aff\,$K=m\le n$. Then

(1) $K\subseteq conv_\eps K$,

(2) $conv_\eps K$ is a body in $\bbr^m={\rm aff}\,K$ of a class $\mathcal{K}^\eps_2$;

(3) if $K\subseteq K'$ and $K'\subset{\rm aff}\,K$ is a body of a class $\mathcal{K}^\eps_2$, then
 $conv_\eps K\subseteq K'$.

\noindent
Hence, the bodies of a class $\mathcal{K}^\eps_2$ are characterized by the property $conv_\eps K=K$.
\end{prop}

\begin{exam}\rm An $\eps$-convex hull of a set $\{A,B\}\subset\bbr^n$ of two different points is either a segment $AB$ if $|AB|<1/\eps$, or again $\{A,B\}$ otherwise.

A set $W=\{A,B,C\}\subset\bbr^n$ of three non-collinear points is contained in a 2-dimensional plane
aff\,$W$ identified with $\bbr^2$. Hence, $conv_\eps W$ is an intersection of all bodies in $\bbr^2$ containing $W$, moreover, $conv_\eps W$ is a part of the triangle $\Delta\,ABC$.
If $\eps>0$ is small, $conv_\eps W$ is obtained by cutting three 2-dimensional discs from $\Delta\,ABC$ in the plane $\bbr^2$, this purely \underline{2-dimensional action} does not depend on the dimension~$n$. As a result we obtain a plane (2-D) figure $conv_\eps\{A,B,C\}\subset\triangle ABC$ ("thin triangle", Fig.~\ref{F-triangle-b}), described in (a), (b), (c) in what follows for some particular cases.

\vskip1mm
 (a) Let $|AB|=|AC|=|BC|=1$. In the case of $\eps\le1$, $conv_\eps W$ is a triangular 2-D domain bounded by three circular arcs of radius $1/\eps$, Fig.~\ref{F-triangle-b}(a). The values $\eps>1$ give us disconnected "$\eps$-convex hulls": $conv_\eps W=W$ when $\eps>\sqrt3$, $conv_\eps W=\{W,O\}$ when $\eps=\sqrt3$,
 and $conv_\eps W$ is a union of $W$ and a small \textit{thin triangle} when $\eps\in(1,\sqrt3)$, Fig.~\ref{F-triangle-b}(b).

\vskip1mm
 (b) Let $1=|AB|>|AC|=|BC|=c>1/2$.
 An outer support to $W$ circle of radius $1/\eps$ through $A,C$ makes angle $\varphi_1=\arcsin(c\,\eps/2)$ with $AC$, and an outer support to $W$ circle of radius $1/\eps$ through $A,B$ makes angle $\varphi_2=\arcsin(\eps/2)$ with $AB$. The requirement $\varphi_1\le\frac12\angle C=\arcsin\frac1{2\,c}$ provides $\eps<1/c^2$.
 From the requirement $\varphi_1+\varphi_2\le\angle A=\arccos\frac1{2\,c}$ it follows
  $$
 \sqrt{(1-c^2\eps^2/4)(1-\eps^2/4)}-c\,\eps^2/4=\cos(\varphi_1+\varphi_2)\ge1/(2\,c),
  $$
 This inequality is equivalent to
 $(1-\frac14c^2\eps^2)(1-\frac14\eps^2)\ge(\frac14c\,\eps^2+\frac1{2\,c})^2$.
 The solution of corresponding equation is $\eps_0=\frac{(4\,c^2-1)^{1/2}}{c(2+c^2)^{1/2}}\le\frac{1}{c^2}$
 for $c\in(\frac12, 1)$.
 Hence $conv_\eps W$ is a thin triangle when $\eps\le\eps_0$ (with zero angles at $A,B$ when $\eps=\eps_0$).
\begin{figure}[ht]
\begin{center}
\includegraphics[scale=.5,angle= 0,clip=true,draft=false]{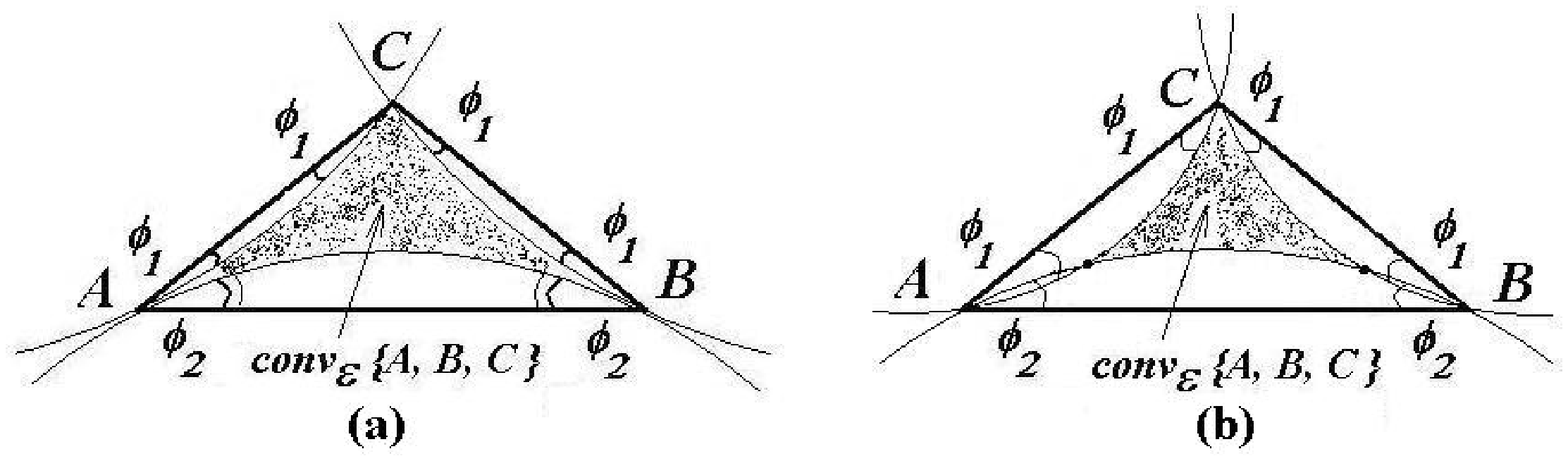}
\caption{\small Thin triangles in $\bbr^2$}
 \label{F-triangle-b}
\end{center}
\end{figure}

\vskip1mm
 (c) Let $c=|AB|<|AC|=|BC|=1$. Similar to (b) we obtain $\varphi_1=\arcsin(\eps/2)$ and $\varphi_2=\arcsin(c\,\eps/2)$. The requirement $\varphi_1\le\frac12\angle C=\arcsin\frac{c}{2}$ provides $\eps<c$.
 From the requirement $\varphi_1+\varphi_2\le\angle A=\arccos\frac{c}{2}$ it follows
$$
 \sqrt{(1-c^2\eps^2/4)(1-\eps^2/4)}-c\,\eps^2/4=\cos(\varphi_1+\varphi_2)\ge c/2,
$$
The solution of corresponding equation is $\eps_0=\frac{(4-c^2)^{1/2}}{(1+2\,c^2)^{1/2}}\ge c$ for $c\in(0, 1)$. Hence $conv_\eps W$ is a thin triangle when $\eps\le c$ (with zero angle at $C$ when $\eps=c$).
\end{exam}

\begin{prop}\label{P-26}
 Let $K\subset\bbr^n$ be a body, and $\eps'>\eps$. Then

 (a) $conv_{\eps'}K\subset conv_\eps K$;\quad
 (b) $K\in\mathcal{K}^\eps_2$ $\ \Longrightarrow K\in\mathcal{K}^{\eps'}_2$.
\end{prop}

\begin{defi}\label{D-eps-convp}\rm
If $n+1$ points $W=\{a_i\}_{i=1}^{n+1}$ in $\bbr^n$ are affine independent and $\eps>0$ is sufficiently small, then $conv_\eps\,W$ is called a \textit{thin $\eps$-simplex}. (Indeed, $conv_\eps\,W$ is homeomorphic to a simplex $conv\,W$). Let $W=\{a_1,\dots a_k\}$ be a finite set of $k>n$ points in $\bbr^n$, and dim\,aff\,$W=m\le n$. Identify aff\,$W=\bbr^m$. Remark that a convex polytope $conv\,W$ in $\bbr^m$ is homeomorphic to a ball, and its vertices form a subset of $W$. If $\eps>0$ is sufficiently small, we call $conv_\eps\,W\subset conv\,W$ a \textit{thin $\eps$-polytope}; it is contained in $\bbr^m$, has the same vertices as $conv\,W$.
\end{defi}

\begin{prop}\label{P-nsimplex} Let $W=\{a_i\}_{i=1}^{n+1}\subset\bbr^n$ be the vertices of a regular simplex $\Delta_W$ with unit edge. If $\eps^2\le\frac{2}{n(n-1)}$ then $conv_\eps W{\subset}\,\Delta_W$ is a thin $\eps$-simplex (a domain with the same vertices bounded by $n+1$ parts of hyper-spheres of radius~$1/\eps$).
\end{prop}

\begin{figure}[ht]
\begin{center}
\includegraphics[scale=.45,angle= 0,clip=true,draft=false]{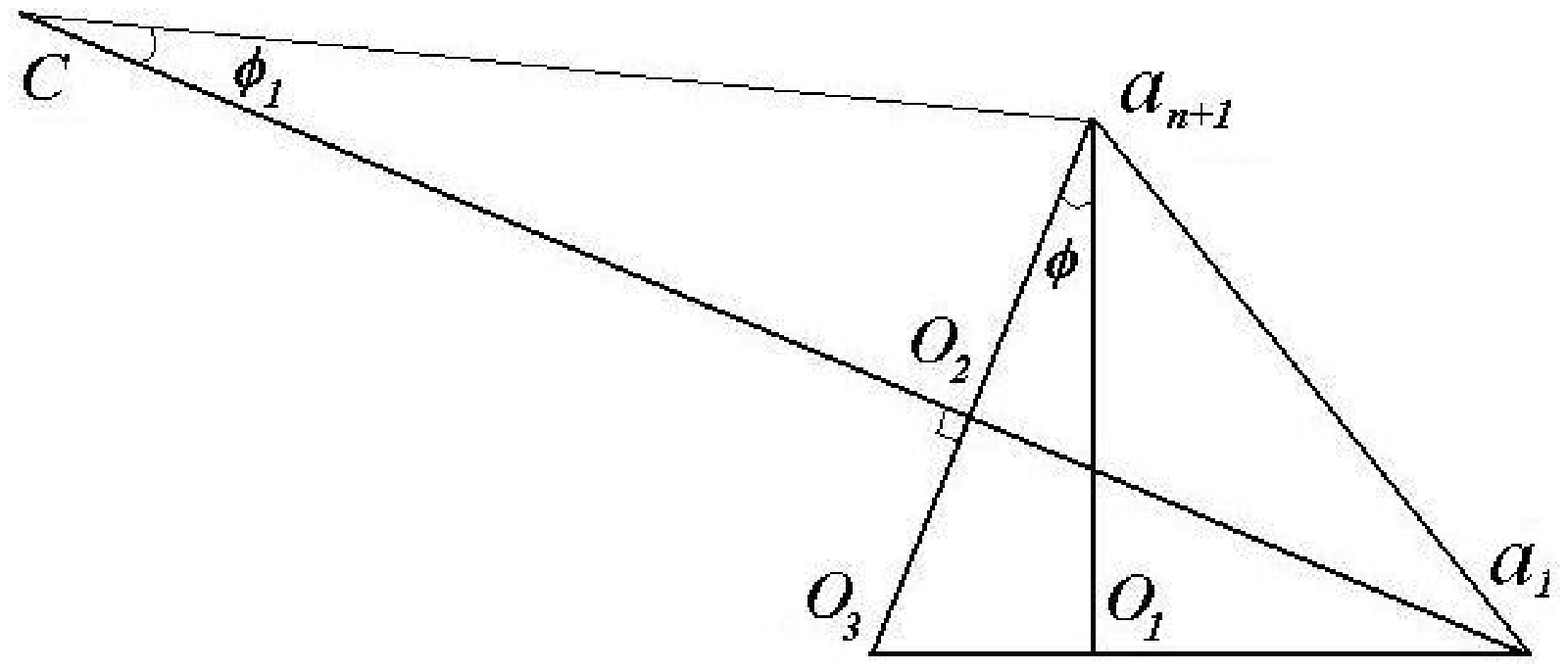}
\caption{\small Estimates for a thin $\eps$-simplex.}
 \label{F-simplex1}
\end{center}
\end{figure}

\textbf{Proof}. Denote by $h_n$ the height of $n$-dimensional simplex with unit edge.
Obviously, $h_2=\sqrt{3}/2$. From the recurrence relation
 $h^2_n+[(n{-}1)/n]^2h^2_{n-1}=1$ by induction it follows the formula $h_n=\sqrt{\frac{n+1}{2n}}$.
 Let $O_1=\frac1n\sum_{i<n+1} a_i$, $O_2=\frac1n\sum_{i>1} a_i$ and $O_3=\frac1{n-1}\sum_{2\le i\le n} a_i$
be the centers of certain faces of~$\Delta_W$, Fig.~\ref{F-simplex1}.
Let $B(C,1/\eps)$ be an outer support ball from the side of the face $\{a_i\}_{2\le i\le n+1}$.
Set $\varphi_1=\angle(O_2 C a_{n+1})$ and $\varphi=\angle(O_2 a_{n+1} O_1)$.
We have $\sin\varphi_1=\frac{n-1}{n}|O_3 a_{n+1}|/|C a_{n+1}|=\frac{n-1}{n}h_{n-1}\eps
 =\eps\sqrt{\frac{n-1}{2n}}$, and hence $\sin\varphi=\frac{|O_1O_3|}{|O_3\,a_{n+1}|}=\frac{1}{n}$. The inequality $\varphi_1\le\varphi$ leads to $\eps^2\le\frac{2}{n(n-1)}$ as required.$\,\square$

\begin{exam}\rm
 Let $W=\{a_i\}_{i\le 4}{\subset}\bbr^3$ be the vertices of a simplex $\Delta_W$.
 According to Proposition~\ref{P-nsimplex}, $conv_\eps W$ is a thin simplex (tetrahedra) for small $\eps>0$. Note that $conv_\eps\{a_1,a_2,a_3\}$ is a thin triangle in the plane aff\,$\{a_1,a_2,a_3\}$, the (relative) interior of $conv_\eps\{a_1,a_2,a_3\}$ is not contained in $conv_\eps W$.
 \end{exam}

\begin{defi}\label{D-Conv}
Let $K\subset\bbr^n$. We call $Conv_\eps K$ the union
 $\bigcup\big\{ conv_\eps\{a_1,\dots a_k\}$: for all $k$ and all sets $\{a_i\}_{i\le k}\in K$ satisfying equality dim\,aff\,$\{a_1,\dots a_k\}={\rm dim\,aff}\,K$\big\}. In other words, $Conv_\eps K$ is a union of $\eps$-convex hulls of all subsets $\{a_1,\dots a_k\}\subset K$ with maximal affine dimension.
\end{defi}

Remark that $Conv_\eps K=\,conv_\eps K$ for a finite set $K\subset\bbr^n$.

\begin{prop}\label{P-Conv}
 $Conv_\eps K=conv_\eps K$ for any compact set $K\subset\bbr^n$ and $\eps>0$.
\end{prop}

\textbf{Proof}. It is enough to prove both inclusions for the bodies $Conv_\eps K,\,conv_\eps K$.

(a) Let a finite set $K'=\{a_1,\dots a_k\}\subseteq K$ has affine dimension $m$, and let a body $K''$ of a class $\mathcal{K}^\eps_2$ contains $K$. Then by Proposition~\ref{P-23}, condition (3), $conv_\eps K'\subseteq K''$. From this and definition it follows $Conv_\eps K\subseteq conv_\eps K$.

(b) Recall that for a compact $K$ and $s>0$ there is a finite set $K^s$ (called \textit{$s$-net of $K$}) such that $\delta(K,K^s)<s$.
Let $K^{(i)}$ be $1/i$-net of $K$, moreover, one may assume $K^{(i)}\subseteq K^{(i+1)},\,\forall i$. The inequality $\delta(K,K^{(i)})<1/i$ means that for any $x\in K$ there is $x_i\in K^{(i)}$ such that $\|x-x_i\|<1/i$. Note that, see Remark~\ref{R-convhull},

-- if $\overset{o}{B^m}(C,1/\eps)\cap K^{(i)}=\emptyset$ for all $i>0$ then
   $\overset{o}{B^m}(C,1/\eps)\cap Conv_\eps K=\emptyset$,

-- if $\overset{o}{B^m}(C,1/\eps)\cap K=\emptyset$ then $\overset{o}{B^m}(C,1/\eps)\cap conv_\eps K=\emptyset$,

\noindent
where as usual, dim\,aff\,$K=m\le n$ and $K\subset\bbr^m$.
We claim that

-- if $\overset{o}{B}(C,1/\eps)\cap K^s=\emptyset$ for all $s>0$
then $\overset{o}B(C,1/\eps)\cap K=\emptyset$.

\noindent
From this and above immediately follows that
 $(\bbr^m\setminus Conv_\eps K)\subseteq(\bbr^m\setminus conv_\eps K)$.
Hence $conv_\eps K\subseteq Conv_\eps K$.

To prove the claim assume an opposite, that there is a ball $int B(C,1/\eps)$ containing a point $x$ from $K$ such that $int B(C,1/\eps)\cap K^{(i)}=\emptyset$ for all $i>0$.
Take natural $N$ such that $s_{i}<1/\eps-\|x-C\|$ for $i>N$.
By definition of $s$-net, for each $i>N$ there is $x_i\in K^{(i)}$ such that $\|x-x_i\|<1/i$. Then $\|x_i-C\|\le\|x_i-x\|+\|x-C\|<1/i+(1/\eps-1/i)=1/\eps$. Hence
$x_i\in\overset{o}B(C,1/\eps)\cap K^{(i)}$ for all $i>N$, a contradiction. $\,\square$

\begin{exam}\label{E-Carat}\rm
The Caratheodory theorem for convex bodies tells us that
{for any compact set $K\subset\bbr^n$ and $x\in conv K$,
there are $m\le n+1$ points $a_1,\dots a_{m}$ in $K$ such that $x\in conv\{a_1,\dots a_{m}\}$}.
Let $W=\{a_i\}_{0\le i\le n}\subset\bbr^{n}$ be the vertices of a simplex with unit edge,
$a'_{0}$ the symmetry of $a_{0}$ relative to aff\,$\{a_1,\dots a_n\}$, and $K=W\cup a'_0$.
Then $O=(a_{0}+a'_{0})/2$ is an inner point of $conv_\eps K$ for small $\eps>0$, but
$O$ does not belong to an $\eps$-convex hull of any $n+1$ points of $K$.
Hence, there is no Caratheodory type theorem for $\eps$-convex bodies of a class $\mathcal{K}^\eps_2$.
\end{exam}

\begin{prob}\rm Find an 'optimal' method for building an $\eps$-convex hull for a finite set $W$ in $\bbr^n\ (n\ge2)$, and estimate its complexity.
\end{prob}

Remark that a most straightforward method, Jarvis's March, which is also
known as {gift wrapping method}, can be easily modified and extended for $\eps$-convex hull of a finite set $W\subset\bbr^2$ with $d_W<2/\eps$.

\subsection{$(\eps,k)$-Convex bodies in normed linear spaces}

We select inside the classes $\mathcal{K}^{\eps}_i$ of bodies in normed linear spaces some smaller classes $\mathcal{K}^{\eps,k}_i$, where $0<k< n$
(for $k=0$ we naturally have $\mathcal{K}^{\eps,0}_i=\mathcal{K}^{\eps}_i$). 
Simply saying, the role of supporting hyperplanes for such bodies will play the cylindrical hyper-surfaces of radius~$1/\eps$.

\vskip1.5mm
Let $P^k\subset\bbr^n$ be a $k$-dimensional plane (i.e., an affine $k$-subspace). Denote~by
\begin{equation}\label{E-tube}
 B(P^k,r)=\{x\in\bbr^n: {\rm dist}(x,P^k)\le r\}
\end{equation}
a {solid cylinder of radius $r>0$ with the axis $P^k$}, and $S(P^k,r)$ its boundary (a cylindrical hypersurface with the same axis). For $k=0$ the formula (\ref{E-tube}) gives us a ball $B^{n}(C,r)\subset\bbr^n$. For $k>0$, $B(P^k,r)$ is a metric product of $\bbr^k$ and a ball $B^{n-k}(C,r)\subset\bbr^{n-k}$. We call $B(P^k,r)$ an \textit{outer support cylinder} of a body $K$ if it doesn't intersect int\,$K$ and intersects $K$ at its boundary points.

\vskip1.5mm
We will extend Definition~\ref{D-eps-conv} of $\eps$-convex bodies.

\begin{defi}\label{D-eps-convk}\rm
A body $K$ in a normed linear space
(for instance $K\subset\bbr^n$) is called $(\eps,k)$-\textit{convex} of a class $\mathcal{K}^{\eps,k}_i$ (for some $\eps>0$) if

\vskip.5mm\hskip-1mm
$\mathcal{K}^{\eps,k}_1$:
{any point $x\in\partial K$ belongs to an outer support cylinder $B(P^k,1/\eps)$ of~$K$},

\vskip.5mm\hskip-1mm
$\mathcal{K}^{\eps,k}_2$:
{any point $x\notin K$ belongs to a cylinder $B(P^k,1/\eps)$ such that}

{int\,$K\cap B(P^k,1/\eps)=\emptyset$},

\vskip.5mm\hskip-1mm
$\mathcal{K}^{\eps,k}_3$:
{any outer support
 ball $B(y,r)$ of $K$ of radius $r<1/\eps$ is contained in}

{an outer support cylinder $B(P^k,1/\eps)$}.

\vskip.5mm\noindent
A~connected boundary component of an $(\eps,k)$-convex convex body will be called an $(\eps,k)$-\textit{convex hypersurface} of a certain class listed above.
\end{defi}

\begin{prop}\label{P-Ki-epsk4}
The class $\mathcal{K}^{\eps,k}_3$ of bodies in complete metric spaces coincides with the class defined by the following weaker condition:

\vskip.5mm
\hskip-1mm
$\mathcal{K}^{\eps,k}_4$: {any
 ball $B(y,r)$ of radius $r<1/\eps$ disjoint from $K$ is contained in}

 {a cylinder $B(P^k,1/\eps)$ that does not intersect $int\,K$}.
\end{prop}

\textbf{Proof} is similar to the proof of Proposition~\ref{P-Ki-eps4}.
$\,\square$

\begin{rem}\rm (a)
One may verify (applying just the set theory arguments and Proposition~\ref{P-Ki-epsk4}) that if a body $K$ is the intersection of (connected) bodies of a class $\mathcal{K}^{\,\eps,k}_2$ (or $\mathcal{K}^{\,\eps,k}_1$) then $K$ also belongs to $\mathcal{K}^{\eps,k}_2$ (resp., $\mathcal{K}^{\,\eps,k}_1$).

(b) Note that $\mathcal{K}^{\,\eps,0}_i=\mathcal{K}^{\,\eps}_i$.
 Indeed, $\mathcal{K}^{\,\eps}_i\subseteq\mathcal{K}^{\,\eps,k}_i$ but the converse is wrong: consider a cube $\{|x_j|\le1,\ 1\le j\le n\}$ and cut off a ball $B(C,1/\eps)$ centered at $C=(\sqrt{1/\eps^2-1},0,\dots 0)$. Replacing cylinders and balls of radius $1/\eps$ in Definitions~\ref{D-eps-convk}, \ref{D-eps-vidk} by half-spaces and half-planes, resp., we obtain (when the bodies are assumed connected) the well-known class $\mathcal{K}^0$ of convex bodies.
\end{rem}

The next definition is inspired by Mazur's theorem (Theorem A.2.1, \cite{nic})
and by the following definition (see Definition 1.1.4 and the certain results in~\cite{vgol}).
A~body $K\subset\bbr^n$ is $(n-1)$-\textit{visible} (or simply, \textit{visible}) if

 $\mathcal{V}$:
 {for any $m$, $ 2 \leq m \leq n$, each $(n-m)$-dimensional plane $P^{n-m}$ which is}

 {disjoint from $K$, belongs to a hyperplane $P^{n-1}$ disjoint from $K$ as well}.

\noindent
Here the zero-dimensional plane is a point.
It is known that a connected $(n-1)$-visible body in $\bbr^n$ is convex.

\begin{defi}\label{D-eps-vidk}\rm
A body $K$ in a $n$-dimensional normed linear space will be called $\eps$-\textit{visible} (of a class $\mathcal{V}^{\,\eps}$ for some $\eps>0$) if

\vskip1.0mm\hskip-1mm
$\mathcal{V}^{\,\eps}$:
 {for any $m$, $ 2 \leq m \leq n$, each plane $Q^{n-m}$ which is disjoint from $K$, belongs}

 {to a cylinder $B(P^{n-m},1/\eps)$
 which is disjoint from int$\,K$}.
\end{defi}

\begin{prop}\label{P-epsk} The following inclusions hold for classes of bodies in $n$-dimensi\-onal Banach spaces:
\begin{equation*}
 \mathcal{K}^{\,\eps,i}_s\subset\mathcal{K}^{\,\eps,j}_s;\quad
 \mathcal{K}^{\,\eps_1,k}_s\subset\mathcal{K}^{\,\eps_2,k}_s,\quad
 \mathcal{V}^{\,\eps_1}\subset\mathcal{V}^{\,\eps_2}\ \
 ({\rm for}\ i>j, \ s=1,2,3, \ \eps_1<\eps_2).
\end{equation*}
\end{prop}

\textbf{Proof} is straightforward. Note that the inclusions 
$ \mathcal{K}^{\,\eps_1}_s\subset\mathcal{K}^{\,\eps_2}_s$ (i.e., $k=0$) 
for $s=1,2,3$ complete Theorem~\ref{T-Ki-eps}.$\,\square$

\begin{theo}\label{T-Ki-epsk} The following inclusions are satisfied for classes of bodies in $n$-dimensional Banach spaces:
\begin{equation*}
 \mathcal{K}^{\,\eps,k}_3\subset\mathcal{K}^{\,\eps,k}_2\subset\mathcal{K}^{\,\eps,k}_1,\qquad
 \mathcal{V}^{\,\eps}\subset\mathcal{K}^{\eps}_2.
\end{equation*}
\end{theo}

\textbf{Proof} is similar to the proof of Theorem~\ref{T-Ki-eps}.
We only show that the inclusions are strong.

(a) Let $K'\in\mathcal{K}^\eps_1\setminus\mathcal{K}^\eps_2$ be a body homeomorphic to a ball in $\bbr^{n-k}$, described in the proof of Theorem~\ref{T-Ki-eps}(a), see Fig.~\ref{F-cond-ia} for $n-k=2$. Then a body $K=K'{\times}[0,1]^k\subset\bbr^{n-k}\times\bbr^k=\bbr^n$ is homeomorphic to a ball, and $K\in\mathcal{K}^{\eps,k}_1\setminus\mathcal{K}^{\eps,k}_2$.

(b) Let $K'\in\mathcal{K}^{\,\eps}_2\setminus\mathcal{K}^{\,\eps}_3$ be a "small" body in $\bbr^{n-k}$, described in the proof of Theorem~\ref{T-Ki-eps}(b), see Fig.~\ref{F-exam2b}. Then a "small" body $K=K'\times[0,1]^k\subset \bbr^{n-k}\times\bbr^{k}=\bbr^n$ is homeomorphic to a ball, and $K$ belongs to $\mathcal{K}^{\eps,k}_2\setminus\mathcal{K}^{\eps,k}_3$.

(c) By definition, $\mathcal{V}^{\,\eps}\subseteq\mathcal{K}^{\,\eps}_2$. The inclusion is
strong, because $(\mathcal{V}^{\,\eps})$ is not valid for $K$ in above (b), and a line
$Q^1\,||\,\bbr^{k}$ through $C'$ (see also Remark~\ref{R-visible}).$\,\square$

\begin{rem}\label{R-visible}\rm
Similar to $\mathcal{K}^{\,\eps,k}_i,\mathcal{V}^{\,\eps}$ classes of bodies may be defined in space forms of non-zero curvature if appropriate cylinders are chosen, and in complex (qua\-ternion) $n$-dimensional normed linear spaces.
\end{rem}

\begin{theo}\label{T-Ki-epsk4}
The class $\mathcal{K}^{\,\eps,k}_3$ of bodies in $\bbr^n$ coincides with the following class:

\vskip.5mm\hskip-1mm
$\mathcal{K}^{\,\eps,k}_{5}$:
 {for every plane $Z^k\in\bbr^n\setminus K$ such that $\rho={\rm dist}(Z^k,K)<1/\eps$,}

 {the intersection $B(Z^k,\rho)\cap K$ is contained in $Z^k+z$ for some $z\in\bbr^n$}.
\end{theo}

\textbf{Proof} is similar to the proof of Theorem~\ref{T-Ki-eps4}.$\,\square$

\begin{prop}
Let a $C^2$-regular hypersurface $\partial K\subset\bbr^n$ bounds a body $K\in\mathcal{K}^{\eps,k}_1$, and $\bn$ be a unit normal to $\partial K$ directed inside. Denote by $\kappa_1(\bn)\le\ldots\le \kappa_{n-1}(\bn)$ the principal curvatures of $\partial K$ with respect to $\bn$.
Then $\kappa_i(\bn)\ge-\eps$ for $i<n-k$ and $\kappa_{n-k}(\bn)\ge0$.
\end{prop}

\begin{prop}
 If $K\in\mathcal{K}^{\,\eps,k}_1\ (k>0)$ is connected then $\partial K$ is connected.
\end{prop}

\textbf{Proof}. Assume a contrary, that $K$ is connected and its boundary $\partial K$ is not. Then a compact component of $\partial K$ bounds a domain $L\subset\bbr^n\setminus K$. The condition $(\mathcal{K}^{\,\eps,k}_1)$ is not satisfied for any $x\in L$ and $k>0$.$\,\square$

\begin{exam}\rm We present a body $K\in\mathcal{K}^{\,\eps,1}_2\ (\eps<1)$ such that $K'=\bbr^4\setminus K$ is not simply connected. In (a) we build a "large" body, and in (b) -- a "small" body.

(a) Let $Q=\{|x_i|\le a,\ i\le 4\}$ be a cube, and $B_1=B(P_{x_1,x_2},b)$, $B_2=B(P_{x_3,x_4},b)$ be cylinders of radius $b\in(1,1/\eps)$ in $\bbr^4$ over coordinate 2-planes. Here $a>b>1/ \eps$.
 Obviously, a "large" body $K=Q\setminus(B_1\cup B_2)$ belongs to a class $\mathcal{K}^{\eps,1}_2$. We claim that $K'=\bbr^4\setminus K$ is not simply connected.
To prove the claim, consider $G=S^3\cup D_{x_1,x_2}\cup D_{x_3,x_4}$, where $D_{x_1,x_2},D_{x_3,x_4}$ are 2-D discs of radius $2\,a$ in corresponding coordinate planes and $S^3$ is a 3-sphere of radius $2\,a$, all with centers at the origin. Remark that $G$ is a homotopic retract of $K'$, because $D_{x_1,x_2}$ is retract of $B_1$, $D_{x_3,x_4}$ is retract of $B_2$, and $S^3$ is retract of $\bbr^4\setminus Q$.
 Now, let $a_1=(2a,0,0,0)$ and $a_4=(0,0,0,2a)$ be two points on the sphere $S^3$. Consider a loop $\gamma\subset K'$ composed by two radii $O a_1$ and $O a_4$ and the short part of meridian of $S^3$ which joins $a_1$ and $a_4$. It is obvious that the loop $\gamma$ is not contractible in $K'$.

(b) Similar example can be constructed for "small" bodies $K$ (i.e., $d_K{<}1/\eps$).
Let $a>1/\eps>c>\sqrt{a^2{-}1}$. Let $Q_1=\{|x_i|\le 1,\ i\le 4\}$ be a cube,
$P_1=\{x_1=c,\,x_3=0\}$ the plane, $C_1=B(P_1,a)$ the cylinder $(x_1{-}c)^2{+}x_3^2\leq a^2$,
$P_2=\{x_1={-}c,\,x_4=0\}$ the plane and $C_2=B(P_2,a)$ the cylinder $(x_1{+}c)^2{+}x_4^2\leq a^2$ in~$\bbr^4$. Denote by $C^2_1$ and $C^2_2$ the intersections of $C_1$ and $C_2$ with $\bbr^3:\ x_2=0$, respectively. These $C^2_1$ and $C^2_2$ cut a hole off the 3-D cube $Q_2=\{|x_i|\le1\ (i=1,3,4); \,x_2=0\}$, and the axis $Ox_1$ traverses $Q^3_2$ through the hole and belongs to the union $C^2_1\cup C^2_2$.
So, it is not contained in $Q^3_2\setminus (C^2_1\cup C^2_2)$. Hence, the plane $Ox_1x_2$ is not contained in the body $Q_1\setminus(C_1\cup C_2)$. This is analogue of cutting the cylinder $C_2$ off the cube in the "large" example (a).

Similarly, let $P_3=\{x_4=c,\,x_2=0\},\,P_4=\{x_4=-c,\,x_1=0\}$ be the planes and $C_3=B(P_3,a):\,(x_4-c)^2+x_2^2\leq a^2,\,C_4=B(P_4,a):\,(x_4+c)^2+x_1^2\leq a^2$ the cylinders in~$\bbr^4$. Denote by $C^3_3,C^3_4$ the intersections of $C_3,C_4$ with $\bbr^3=\{x_3=0\}$, respectively. These $C^3_3,C^3_4$ cut off a hole in the 3-D cube $Q_3=\{-1\leq x_1,x_2,x_4 \leq 1,\ x_3=0\}$, the axis $Ox_4$ traverses $Q_3$ through the hole and is contained in $C^3_1\cup C^3_2$. So, it is not contained in $Q_3\setminus(C^3_3\cup C^3_4)$.
Hence, the plane $Ox_3x_4$ does not belong to the body $Q_1\setminus(C_3\cup C_4)$. This is analogue of cutting the cylinder $C_1$ off the cube in~(a).
 If we cut all these four "large" cylinders from a "small" cube $Q_1$, we obtain a body $K_1$ homeomorphic to "large" $K$ described in the part (a). So, $\bbr^4\setminus K_1$ is not simply-connected.
\end{exam}

\section{Applications to Geometric Tomography}
\label{sec:3}

We introduce the circular projections (Section~\ref{sec:circ}), and apply them to the problem of determination of $\eps$-convex bodies by their projecti\-on-type images (Section~\ref{sec:tomogr}), the results are related to geometric tomography due to \cite{rgar06}.

\subsection{Circular projections}\label{sec:circ}

The orthogonal projections of different $\eps$-convex bodies onto all hyperplanes may coincide (as for a cube and $\eps$-convex hull of its 12 edges for small~$\eps$). Hence we need more complicated maps for the role of projections of such bodies.



\begin{defi}\label{D-circ-pro}\rm
 Given $\omega,C\in\bbr^n$ denote by $l(C,\omega)$ the straight line through $C$ in direction $\omega$,
 and $P(C,\omega)$ the hyperplane through $C$ and orthogonal to~$\omega$.
 We define a map $f_{C,\,\omega}:(\bbr^n\setminus l(C,\omega))\to(P(C,\omega)\setminus C)$
 called a \textit{circular projection} onto $P(C,V)$, as follows.
 Take any $x\in\bbr^n\setminus l(C,\omega)$.
 Let $S^1_x$ be a circle through $x$ and centered at $C$, whose plane is parallel to vectors $\overrightarrow{Cx}$ and $\omega$. Then $f_{C,\,\omega}(x)$ is the nearest to $x$ point of intersection $S^1_x\cap P(C,\omega)$.  Remark that $f_{C,\,\omega}(x)=x$ if $x\in P(C,\omega)\setminus C$.
 \end{defi}

The circular projection $f_{C,\,\omega}$  can be represented explicitly by a formula.
If $C=O$ then
\begin{equation}\label{E-fomega}
\begin{array}{c}
 f_{O,\,\omega}(x)=\frac{\|x\|}{\sqrt{x^2-(\omega,x)^2}}[x-(\omega,x)\,\omega],\quad
 \forall x\hskip-1pt\not\in\,l(O,\omega).
\end{array}
\end{equation}
A circular projection $f_{C,\,\omega}$ is a smooth map, it keeps the distance from points to $C$ and maps bodies (not intersecting $l(C,\omega)$) into bodies of $P(C,\omega)$.
One may easily extend Definition~\ref{D-circ-pro} using non-planar "screens" $P(C,\omega)$ of the projections.

\begin{rem}\rm
The circular projections are related to hyperbolic geometry modeled on a half-space
(Poincar\'{e} model).
Denote by $l^+(C,\omega)$ a ray defined by $C$ and a direction $\omega$,
$P^+(C,\omega)$ an open half-space bounded by a plane $P(C,\omega)$ and containing $l^+(C,\omega)$. Consider the standard \textit{hyperbolic metric} in $P^+(C,\omega)$.
Then trajectories of $f^+_{C,\,V}:P^+(C,\omega)\to(P^+(C,\omega)\setminus C)$
(the restriction of a circular projection) become geodesics orthogonal to the fixed geodesic
$l^+(C,\omega)$.
\end{rem}

In the next definition we are based on the notion of apparent contour of a surface under orthogonal projection onto the plane, see \cite{vgol}, \cite{haef}.

\begin{defi}\label{R-apparentf}\rm
Let $K\subset\bbr^n$ be a compact body bounded by a smooth hypersurface $\partial K$.
Let $f=f_{C,\omega}$ be a circular projection onto a hyperplane $P(C,\omega)$
such that $K\cap\,l(C,\omega)=\emptyset$.
The set of points $y\in P(C,\omega)$ such that a circle $S^1_y$ (i.e. a trajectory of $f$ containing $y$) is tangent to $\partial K$ at some point $z\in\partial K$ is called the \textit{apparent contour} of the hypersurface $\partial K$ under a circular projection $f$ and is denoted by $C(\partial K,f)$.
 \end{defi}

\begin{exam}\label{Ex-ball}\rm
For $n=3$ consider a circular projection $f=f_{C,\omega}$ onto coordinate plane $P(C,\omega)=\{z=0\}$, where $C=O$, $\omega=(0,0,1)$, see Fig.~\ref{F-ball1}.
We obtain
\begin{equation}\label{E-fomega3}
 f(x,y,z)=(\alpha\,x,\alpha\,y,0),
 \quad{\rm where}\quad
 \alpha=\sqrt{1+z^2/(x^2+y^2)}.
\end{equation}
\begin{figure}[ht]
\begin{center}
\includegraphics[scale=.37,angle= 0,clip=true,draft=false]{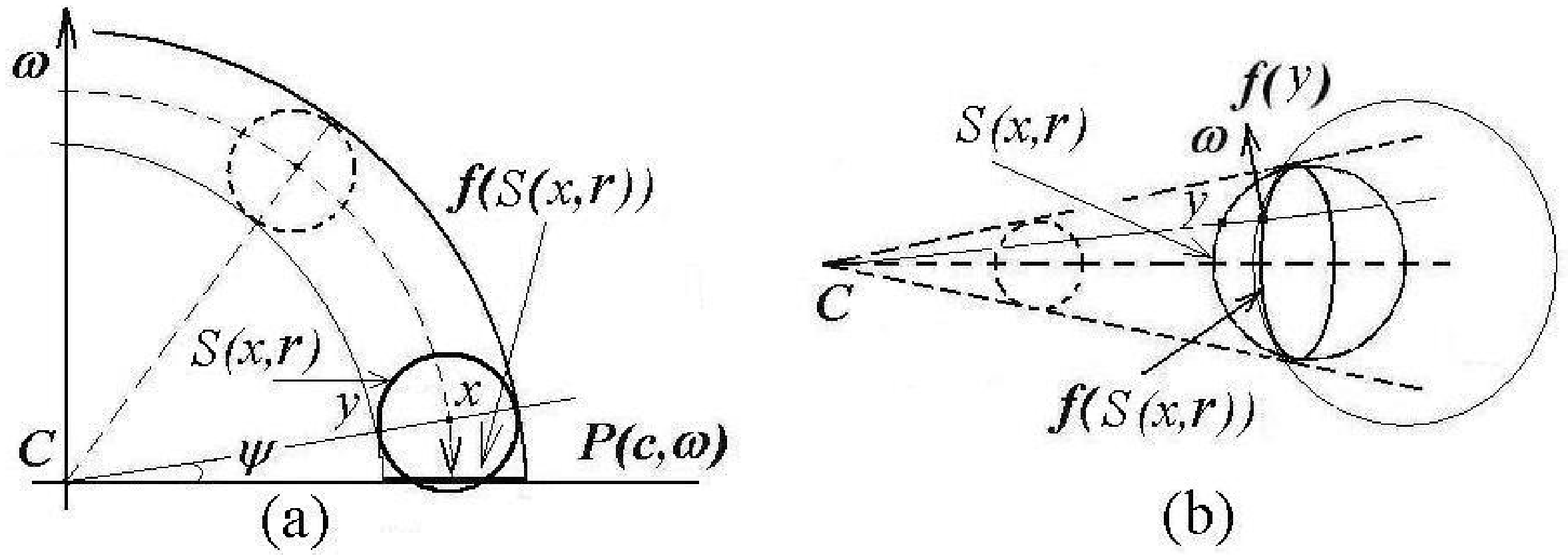}
\caption{\small Image of a sphere under a circular projection $f$: side and upper views.}
 \label{F-ball1}
\end{center}
\end{figure}
Let $S(A,1)\subset\bbr^3$ be a unit sphere with center at $A(0,2,0)$. A circle $S^1(A,1)=S(A,1)\cap P(C,\omega)$ has the following parametrization: $\gamma_0: \br_0(t)=[\sin\,t,\,2+\cos\,t,\,0]$.
The image of $S(A,1)$ under rotation $R_x(\psi)$ about the $x$-axis by an angle $\psi\in(0,\pi/2]$ is a sphere $S(A_\psi,1)$ centered at $A_\psi=R_x(\psi)(A)$.
A great circle $S^1(A_\psi,1)=R_x(\psi)(S^1(A,1))$ on $S(A_\psi,1)$ projects onto the apparent contour $\gamma_\psi$ of $S(A_\psi,1)$, hence $f(S(A_\psi,1))$ is bounded by the curve $f(S^1(A_\psi,1))$.
We~will use (\ref{E-fomega3}) to compute the curvature of the images $f(B(A_\psi,1))$, see Figs.~\ref{F-circ_ball}(a-c) for eight values $\psi\in\{0.15, 0.3, 0.45, 0.524, 0.6, 0.75, 0.9, 1.05\}$.
Remark that $\pi/6\approx0.524$. The parametrization of $\gamma_\psi$ is the following:
\begin{equation}\label{E-gamma3}
\begin{array}{c}
 \br_\psi(t)=\sqrt{\frac{(2+\cos\,t)^{2}+\sin^2t}{\sin^2\psi(2+\cos\,t)^{2}+\sin^2t}}
 \,\big[\sin\,t,\,(2{+}\cos\,t)\cos\,\psi\big].
\end{array}
\end{equation}
\begin{figure}[ht]
\begin{center}
\includegraphics[scale=.34,angle= 0,clip=true,draft=false]{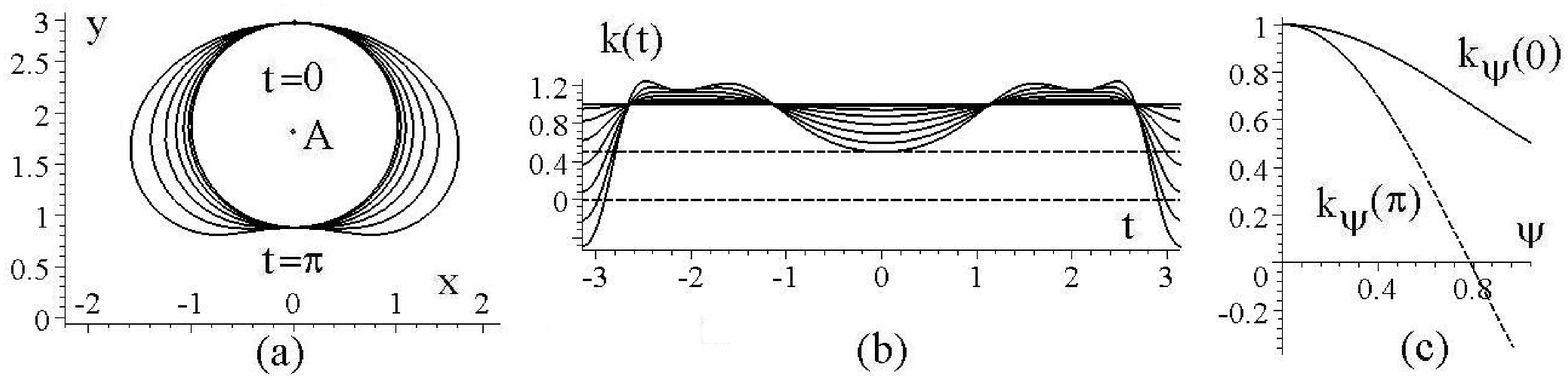}
\caption{\small Example~\ref{Ex-ball}: (a) Image of a sphere $S_\psi(A,1)$ under a circular projection.
(b)~Curvature of $\gamma_\psi: t\to\br_\psi(t)$.
(c)~Comparison of curvature at $\br_\psi(0)$ and $\br_\psi(\pi)$.}
 \label{F-circ_ball}
\end{center}
\end{figure}
The curvature of $\gamma_\psi: t\to\br_\psi(t)$ is positive for $\psi\le\psi_0\approx0.75$, hence the curve itself is convex for these values of $\psi$, see also Figs.~\ref{F-circ_ball}(a,b) for the values $\psi\in\{0.15, 0.3, 0.45, 0.524, 0.6, 0.75\}$. One may conclude that a domain $f(S_\psi(A,1))$ is convex for $\psi\le\psi_0$.
The curvature of $\gamma_\psi$ is greater than 0.5 for $\phi\le\pi/6$, see also Fig.~\ref{F-circ_ball}(b) for the values $\psi\in\{0.15, 0.3, 0.45, 0.524\}$. Hence  $\gamma_\psi$ can be rolled in $\bbr^2$ without sliding inside a disc of radius $2$ (see \cite{top}, Problem 1.7.10). It is easy to verify that the curvature of $\gamma_\psi: t\to\br_\psi(t)$ takes its minimum at $t=\pi$, see Figs.~\ref{F-circ_ball}(b,c).
\end{exam}

\begin{defi}\rm
Consider a curve $\gamma$ that is either contained in the interior of a half-space $\bbr^2_+\subset\bbr^{n}$ or meets its boundary $m$ orthogonally. Now, the group $SO(n)$ contains a subgroup $G$ isomorphic to $SO(n{-}1)$ which acts on the hyperplane orthogonal to $m$. This gives rise to a $G$-invariant hypersurface $M^{n-1}\subset\bbr^n$ called a \textit{hypersurface of revolution about $m$ with $\gamma$ as profile}.
\end{defi}

The group $G$ acts on $M^{n-1}$ and the orbits of this action are just $(n-2)$-dimensional spheres.
In particular, they are umbilical submanifolds of $M^{n-1}$, i.e., at any point their principal curvatures coincide, see, for example, \cite{toj_06}. Remark that $m$ is the symmetry axis of a curve $M^1$ when $n=2$, and $m$ is the usual rotation axis of a surface $M^2$ when $n=3$.
Hence, there are two distinct principal curvatures at any point on a hypersurface of revolution $M^{n-1}\subset\bbr^n\ (n>2)$: one of multiplicity 1 is the curvature of a profile curve (i.e., the boundary curve of intersection with a half-plane determined by $m$ and the point); another one has multiplicity $n-2$.

\begin{lem}\label{L-Top}
Let $M\subset\bbr^n$ be a smooth hypersurface of revolution about axis $m$ with a profile curve $\gamma$, whose curvature at each point is not smaller than $1/a$. Then $M$ can be rolled without sliding inside a ball of radius $a$.
\end{lem}

\textbf{Proof}. For any $P\in M$ consider profile $\gamma_P=\bbr^2_+\cap M$ through $P$.
Then $\tilde\gamma=\bbr^2\cap M$ is a smooth closed curve in $\bbr^2$ (symmetric with respect to $m$) whose curvature at each point is not smaller than $1/a$. Hence, $\tilde\gamma$ can be rolled in $\bbr^2$ without sliding inside a disc of radius $a$, see \cite{top}.
We will prove that a surface $M$ has a similar property relative to balls of radius $a$.

1. Locate a circle $S^1(x,a)\subset\bbr^2$ of radius $a$ so that $m$ is the $y$-axis and the origin $O$ belongs to $S^1(x,a)$ and $\tilde\gamma$. Since $M$ is smooth, the tangent line to $\tilde\gamma$ at $O$ coincides with the $x$-axis. Then $\tilde\gamma$ is contained in $S^1(x,a)$.

Introduce the arc length parameter $s$ counted from $O$ on $\tilde\gamma$ and $S^1(x,a)$.
Then $k(s)\ge a$ and $\alpha(s)=\int_0^s k(s)\,ds\ge\frac sa$. The equations of $\tilde\gamma$ and $S^1(x,a)$ are
\begin{equation*}
  \tilde\gamma:
  \begin{array}{c}
  x_1=\ \int_0^s \cos\alpha(s)\,ds,\\
  y_1=\int_0^s \sin\alpha(s)\,ds,\\
  \end{array}\quad
  S^1(x,a):
  \begin{array}{c}
  x_2=\ \int_0^s \cos(s/a),ds,\\
  y_2=\int_0^s \sin(s/a)\,ds.\\
  \end{array}
\end{equation*}
Let at the points $P=P(s_1)\in\gamma_P$ and $Q(s_2)\in S^1(x,a)$ the tangent lines are parallel (make an angle $\theta$ with the $y$-axis), and thus $s_2/a=\alpha(s_1)$. We obtain
\begin{equation*}
\begin{array}{c}
 x_1(s_1)=\int_0^{s_1}\cos\alpha(s)\,ds\le\int_0^{s_1}\cos\frac sa\,ds=a\sin\frac{s_1}a\le a\sin\alpha(s_1).
\end{array}
\end{equation*}
On the other hand,
 $x_2(s_2)=\int_0^{s_2}\cos\frac sa\,ds=[\,a\sin\frac sa\,]\,|_0^{a\alpha(s_1)}=a\sin\alpha(s_1)$.
It~is shown that $x_1(s_1)\le x_2(s_2)$.
By Meusnier's theorem the principal curvature at $P$ of a surface of revolution $M$ along a parallel equals to $\cos\theta/x_1(s_1)$ that is greater than $\cos\theta/x_2(s_2)=1/a$, i.e., the normal curvature of a sphere $S(x,a)$.

2. Consider a sphere $S(y,a)$ that has a common inner normal with $M$ at $P\in\gamma_P$.
Then $y\in\bbr^2$. Let a circle $S^1(y,a)$ be an intersection $S(y,a)\cap\bbr^2$.
As was shown above, a point $P\in\gamma_P$ is less distanced from the axis $m$ than a point
$Q\in S^1(x,a)$ at which the tangent direction is the same as of $\gamma_P$ at $P$. Hence $y$ does not belong to the half-plane $\bbr^2_+$ containing $\gamma_P$ and bounded by the axis $m$. In other words, $\bbr^2_+$ includes an arc $\tilde S$ of $S^1(y,a)$ smaller than a semi-circle. The profile $\gamma_P$ is contained in the domain of $\bbr^2_+$ bounded by $\tilde S$, hence $\tilde M$, the surface of revolution of $\tilde S$ about $m$, contains $M$. Obviously, $\tilde M\subset B(y,a)$. Hence, a ball $B(y,a)$ includes $M$.$\,\square$

\begin{lem}\label{L-Aball1}
Let $f=f_{C,\,\omega}$ be a circular projection onto a hyperplane $P(C,\omega)\subset\bbr^{n}$. If a sphere $S(x,r)$ does not intersect the line $l(C,\,\omega)$, then the boundary $M$ of
$f(S(x,r))$ is a smooth $(n-2)$-dimensional hypersurface of revolution about axis $m_{C,x}$ defined as the line in $P(C,\omega)$ through the points $C$ and $f(x)$.
\end{lem}

\textbf{Proof}. Let $D_C\,||\,Cx$ be diameter of the sphere $S(x,R)$.
Given $z\in D_C$ consider $(n-2)$-dimensional plane $V_z$ orthogonal to the vectors $\omega$ and $\overrightarrow{Cx}$ such that $z\in V_z$. For any interior point $z$ of this diameter $D_C$ the circular projection $f$ of $(n-3)$-dimensional sphere $S_z=S(x,r)\cap V_z$ is a smooth surface $S'_z\subset P(C,\omega)$ symmetric with respect to rotations in $P(C,\omega)$ about the axis $m_{C,x}$.
Let $n$ be a unit vector orthogonal to $Cx$ and $V_z$.
The union $\cup_{z\in D_C}S_z$ is a great hypersphere $S_y(x,r){\subset} S(x,r)$ through $y=Cx\cap S(x,r)$
and orthogonal to $n$.
The union $\cup_{z\in D_C}S'_z$ composes the boundary hypersurface $M$ of $f(S(x,r))$.

One may verify this by the following derivations based on (\ref{E-fomega}).
The sphere $S_z$ has radius $r_z{=}\sqrt{r^2{-}|zx|^2}\le r$. Let $X\in S_z$ and $X=z+r_zv$, where $v\in V_z$ is a unit vector. Recall that $(v,\omega)=0$. Assuming $C=O$, we deduce from (\ref{E-fomega}) that
$f(X)=\alpha(z'+r_zv)$, where $z'=z-(z,\omega)\,\omega$. Since $X^2=z^2+r_z^2$, we have that $\alpha^2=\frac{\|X\|^2}{X^2-(\omega,X)^2}=\frac{z^2+r_z^2}{z^2+r_z^2-(\omega,z)^2}$
does not depend on~$v$. Hence $S'_z=S(\alpha\,z',\alpha\,r_z)\cap V_z$ is a $(n{-}3)$-dimensional sphere in $P(C,\omega)$ of radius $\br_1=\alpha\,r_z$ centered at~$\alpha\,z'\in m_{C,x}$.
Moreover, ${\rm aff}(S'_z)\perp m_{C,x}$.
The hypersurface of revolution $M$ is smooth, because the vector function $\br_1(z)\ (0\le z\le r)$ belongs to a class $C^\infty(0,r)$, vanishes at the ends of its domain, $\br_1(0)=\br_1(r)=0$, and its derivatives also vanish, $\frac{d}{dz}\br_1(0_+)=\frac{d}{dz}\br_1(r_-)=0$.$\,\square$

\vskip1.5mm
From Lemmas \ref{L-Top}, \ref{L-Aball1} and derivations of Example~\ref{Ex-ball} we obtain the following.

\begin{lem}\label{L-Aball2}
Let $f=f_{C,\,\omega},\,S(x,r)$ are as in Lemma~\ref{L-Aball1}.
Suppose that $r/|Cx|<1/2$ and the angle $\psi=\angle(Cx,P(C,\omega))<\pi/6$. Then the boundary of $f(S(x,r))$ is a smooth convex hypersurface of revolution $M$ in $P(C,\omega)$ with a profile curve $\gamma$, whose curvature at each point is not smaller than $1/(2r)$. For each point $y\in M$ there is a ball of radius $2\,r$ containing $M$ and with the same tangent plane at $y$.
\end{lem}

\textbf{Proof}. By Lemma \ref{L-Aball1}, $M$ is a smooth hypersurface of revolution about $m_{C,x}$ which is defined as the line in $P(C,\omega)$ through the points $C,\,f(x)$. By~Lemma~\ref{L-Top}, we only need to show that the curvature of profile $\gamma$ is not smaller than $1/(2r)$. By conditions we have $\sin\psi<1/2$. Recall that $S_y(x,r)\subset S(x,r)$ is a great hypersphere through $y=Cx\cap S(x,r)$ and is orthogonal to $n$ (see Lemma~\ref{L-Aball1}).
Similarly to Example~\ref{Ex-ball} (see also Fig.~\ref{F-ball1}), we assume $C=O$, then parameterize any great circle $\gamma_{y}\subset S_y(x,r)$ through $y$, and compute the curvature of its image $\gamma(t)=f(\gamma_{y})$.
 The parametrization of $\gamma_{y}$ is $\br(t)=(\frac{r}{\|x\|}\cos\,t+1)x+r\sin\,t\,v$, where $v$ is a unit vector with the property $v\perp x\wedge \omega$. Using (\ref{E-fomega}) we obtain
\begin{equation*}
 f(\gamma_{y}):\ \br_1(t)=\alpha(t)[(1+\lambda\cos\,t)x'+r\sin\,t\,v],
\end{equation*}
where
$\alpha^2(t)=\frac{\|x\|^2+r^2+2r\|x\|\cos\,t}{\|x\|^2+r^2+2r\|x\|\cos\,t-(\omega,x)^2(\frac{r}{\|x\|}\cos\,t+1)^2}$
and $x'=x-(x,\omega)\,\omega$. Set $\lambda=r/\|x\|$. Then $\alpha^2(t)=\frac{1+\lambda^2+2\lambda\cos\,t}{1+\lambda^2+2\lambda\cos\,t-\sin\psi^2(\lambda\cos\,t+1)^2}$.
Recall that $\|x'\|=\|x\|\cos\psi>\sqrt3\,r$ and $\lambda<1/2$.
In $\bbr^2$ with the orthonormal basis $\{v,x'/\|x'\|\}$ we have
 $f(\gamma_{y}):\,\br_1(t)=r\,\alpha(t)[2(1{+}\lambda\cos\,t)\cos\psi, \sin\,t]$.
 Compare it with the curve of smaller minimal curvature (i.e., $\lambda=1/2,\,\psi=\pi/6$)
\begin{equation*}
\begin{array}{c}
 \br_2(t)=\frac{r\,(5+4\,\cos\,t)^{1/2}}{[5+4\,\cos\,t-(0.5\cos\,t+1)^2]^{1/2}}[(1+0.5\cos\,t)\sqrt3, \sin\,t].
\end{array}
\end{equation*}
Direct calculations show that the curvature of $\frac1r\br_2(t)$ is not less than $1/2$. Hence, the curvature of $f(\gamma_y)$ is not less than $1/(2\,r)$.$\,\square$

\subsection{Determination of $\eps$-convex bodies by projection images}
\label{sec:tomogr}

\begin{defi}\label{D-Eomega}\rm
Let $\Omega\subset S^{n-1}$ be a set of unit vectors that has nonempty intersection with any great $(n-2)$-dimensional sphere.
Given bodies $K,L\subset\bbr^n$ and $\eps>0$, denote by $\mathcal{E}_{\Omega, K, L, \eps}$ the collection of all punctured hyperplanes $P(C,\omega)$ in $\bbr^n$ such that the normal $\omega\in\Omega$ and
 either dist$(C,K)<2/\eps$ and $P(C,\omega)$ intersects $L$
 or dist$(C,L)\le2/\eps$ and $P(C,\omega)$ intersects $K$.

Remark that $C$ is not uniquely determined by a plane $P(C,\omega)$.
 \end{defi}

Next theorem generalizes Lemma 1.2.1 in \cite{vgol}.

\begin{theo}\label{T-0eps}
 Let the bodies $K,L\subset\bbr^n\ (n>2)$ belong to a class $\mathcal{K}^\eps_2$, and $\delta(K,L)<1/\eps$. If for all hyperplanes $P(C,\omega)\in\mathcal{E}_{\Omega, K, L, \eps}$ the images of corresponding circular projections $f_{C,\omega}(K)$ and $f_{C,\omega}(L)$ coincide, then these bodies $K$ and $L$ coincide themselves in the ambient space $\bbr^n$.
\end{theo}

\textbf{Proof}. Suppose an opposite that ${\rm int}\,L\ne{\rm int}\,K$.
Then there is a ball $B(y,r)\in{\rm int}\,L\setminus K$
(modulo change of names $K,L$).
By~condition ($\mathcal{K}^\eps_2$) there is a ball $B(C,1/\eps)$ containing $y$ and not intersecting int\,$K$.

Let us take a point $x\in B(y,r)$ with the property $\rho(x,C)<1/\eps$.
Due to (\ref{E-hausd2}) we have $x\in K_{\delta}$, where $\delta=\delta(K,L)$.
By the triangle inequality, from $\rho(x,C)<1/\eps$ and dist$(x,K)\le\delta<1/\eps$,
it follows that dist$(C,K)<2/\eps$. Moreover, $x\ne C$ (otherwise dist$(C,K)<1/\eps$, hence $B(C,1/\eps)$ intersects int$K$).
 By definition of $\Omega$, there is $\omega\in\Omega$ orthogonal to a nonzero vector $\overrightarrow{xC}$.
Notice that $P(C,\omega)\cap L\not=\emptyset$ because of $x\in P(C,\omega)$.
Hence a hyperplane $P(C,\omega)\in\mathcal{E}_{\Omega, K, L, \eps}$.
Let $f=f_{C,\,\omega}$ be the circular projection onto $P(C,\omega)\setminus C$.

Consider a ball $B^{n-1}=B(C,1/\eps)\cap P(C,\omega)$.
We see that $f(x)=x\in{\rm int}B^{n-1}$ and the image $f(K)$ is contained in the complement of $B^{n-1}$ in $P(C,\omega)\equiv\bbr^{n-1}$.
Hence a point $x\in f(L)$ does not belong to the image $f(K)$ in $P(C,\omega)$,
Fig~\ref{F-Th-4}(a), a contradiction.$\,\square$

\vskip1.5mm
Denote by $\pi:\bbr^n\to P$ the orthogonal projection onto a hyperplane $P$.
In~the sequel we extend and study the following stability problem (\cite{vgol}, \cite{gr87}, \cite{gr97}): {if~for some convex bodies $K,L\subset\bbr^n$, $\eps\ge0$ and every hyperplane $P$ it is known that $\delta_t(f(L),f(K))\le\eps$, what can be said about the size of $\delta_t(L,K)$}?

 \begin{defi}\rm
 If bodies $K,L$ are translates of each other we write $K\simeq L$.
 The~\textit{translative Hausdorff distance} between compact sets $K,L\subset\bbr^n$ is defined as, \cite{gr87}, \cite{gr97},
 \begin{equation*}
 \delta_t(K,L)=\inf\{\delta(K+p,L):\ p\in\bbr^n\}.
 \end{equation*}
  Denote by $Rot_n(C)$ the set of all rotations in $\bbr^n$ about $(n-2)$-dimensional subspaces through a point~$C$. If two bodies $K,L$ are related each to other by a rotation $\phi\in Rot_n(C)$,
 i.e., $L=\phi(K)$, we write $K\overset{C}{\simeq} L$.
 The~\textit{rotational relative to $C$ Hausdorff distance} between compact bodies $K,L\subset\bbr^n$ will be called
 \begin{equation*}
 \delta_{rC}(K,L)=\inf\{\delta(\phi(K),L):\ \phi\in Rot_n(C)\}.
 \end{equation*}
 \end{defi}

 A~\textit{ball $B\subset\bbr^n$ supports} $K$ in the direction $\omega$ if $K\subseteq B\subset H_K^+(\omega)$, \cite{gr97}.
 The~\textit{smallest support ball} in the direction $\omega$ is denoted $B_K(\omega)$,
 its radius is denoted $R_K(\omega)$.
 Clearly if $R_K(\omega)<\infty$, then $\omega$ is a regular direction of~$K$.

\begin{defi}\label{D-Eomega2}\rm
 Given nonzero vectors $\omega_0,\omega_1\in\bbr^n$, we define a set of unit vectors $\Omega_{\omega_0,\omega_1}\subset S^{n-1}$ by the condition $\Omega_{\omega_0,\omega_1}=span\{\omega_1,\omega_0^\perp\cap\,\omega_1^\perp\}$.

 Let $K,L\subset\bbr^n$ be compact bodies.
 Suppose that $R_K(\omega_0)<\infty$ and $A=S_K(\omega_0)\in S_L(\omega_0)$ for some direction $\omega_0$.
 Denote by $\mathcal{E}_{\omega_0, A, K, L, \eps}$ the collection of all punctured hyperplanes $P(C,\omega)$
 intersecting $K$ and $L$ such that
 the normal $\omega\,{\in}\,\Omega_{\omega_0,\overrightarrow{CA}}$
 and $B(C,1/\eps)$ is an outer support ball of either $K$ or~$L$.

\vskip.1mm
Remark that $C$ is not uniquely determined by a plane $P(C,\omega)$.
 \end{defi}

By Definitions~\ref{R-apparentf},\,\ref{D-Eomega2}, the image of a point $A$ under a circular projection~$f_{C,\omega}$ belongs to the apparent contours of a surface $\partial K$ and a sphere $\partial B_K(\omega_0)$. For $\eps\to0$ the family $\mathcal{E}_{\omega_0, A, K, L, \eps}$ reduces to the collection of planes parallel to~$\omega_0$.

Theorem \ref{T-0eps} tells us that if circular projections of $\eps_0 $-convex bodies $K,L$ coincide (for large dist$(C,K)$, dist$(C,L)$, say $>1/\eps_0$), then $K=L$.
The next theorem shows that if these projections are "translation equivalent" with precision $\eps>0$, then $K,L$ should be "translation equivalent" with corresponding precision,
moreover, given $\eps>0$ one should take sufficiently small $\eps_0$ in order to get such an estimate. Similar theorem for convex bodies has been proven in~\cite{gr97}.
 "\dots It is not necessary to consider projections onto all planes of $\bbr^3$ but only onto planes that contain a given line, say $l$, and one additional plane that is orthogonal to $l$.
In various practical situations regarding the determination of bodies from the 'pictures of their shadows', and also from a purely geometric point of view, it is of interest to study if a body can be determined without the knowledge of this exceptional projection onto a plane orthogonal to $l$. In general it is not possible. Consider, for example, two right cylinders of equal height; one having as base a circular unit disc, the other a Reuleaux triangle of width $2$. Then these cylinders are obviously not translates of each other, but have translation equivalent rectangles as 'lateral' projections", \cite{gr97}.

For sufficiently small $\eps>0$ the circular projections of these cylinders are almost translation or rotation equivalent, but the original bodies not. Hence the condition for the radius of $R_K(\omega_0)$ in Theorem~\ref{T-01} is not superabundant.

\begin{theo}\label{T-01}
Let the bodies $K,L\subset\bbr^n\ (n>2)$ belong to a class $\mathcal{K}^{\eps_0}_{3}$ and let $d_K,d_L,\delta(K,L)$ are less than $1/(2\,\eps_0)$. Suppose that $R_K(\omega_0)<1/(3\,\eps_0)$ and $A=S_K(\omega_0)\in S_L(\omega_0)$ for some direction $\omega_0$. If for each plane $P(C,\omega)\in\mathcal{E}_{\omega_0, A, K, L, \eps_0}$ (where
$C$ is selected in (a) of Definition~\ref{D-Eomega2}) the images of a circular projection $f_{C,\,\omega}$ satisfy any of the inequalities
\begin{equation}\label{E-S2}
 a)~\delta_t(f_{C,\,\omega}(K), f_{C,\,\omega}(L))\le\eps,\qquad
 b)~\delta_{rC}(f_{C,\,\omega}(K), f_{C,\,\omega}(L))\le\eps
\end{equation}
for some $\eps\ge0$, then the corresponding inequality holds
 \begin{equation}\label{E-S3}
 \begin{array}{c}
 \hskip-4mm a)~\delta_{t}(K,L)\le 2(\sqrt{2\,R_K(\omega_0)}{+}\sqrt{\eps})\sqrt{\eps},\\
 b)~\delta_{rC}(K,L)\le(\sqrt{8\/R_K(C)}{+}3\sqrt{\eps}{+}\tilde\eps)\sqrt{\eps},
  \end{array}
\end{equation}
where $\tilde\eps$ is given in~Lemma \ref{L-01b} in what follows.
\end{theo}

\textbf{Proof}. a) We first show using condition $d_L,d_K,\delta(K,L)<1/(2\eps_0)$ that
\begin{equation}\label{E-S16}
 \delta(K,L)\le \delta(f_{C,\,\omega}(K), f_{C,\,\omega}(L))
\end{equation}
\underline{for some} circular projection $f_{C,\,\omega}$ onto certain hyperplane $P(C,\omega){\in}\,\mathcal{E}_{\omega_0, A, K, L, \eps_0}$.
To prove this, consider a pair $a\in K,\, b\in L$ such that $\delta=\delta(K,L)={\rm dist}(b,K)=\|a-b\|$ (modulo change of names $K,L$, see below Proposition~\ref{P-ab}). Obviously, $a\in B(b,\delta)$.
By~($\mathcal{K}^{\eps_0}_3$) for any small $\alpha>0$ there is an outer support (of $K$) ball $B(C',1/\eps_0)$ separating $K$ from a ball $B(b,\delta{-}\alpha)$.
We have $C'\to C$ when $\alpha\to0$, hence an outer support (of $K$) ball $B(C,1/\eps_0)$ contains $a$ and separates int$K$ from a ball $B(b,\delta)$. Moreover, $b\in[a,C]$, because the points $a,b$ belong to the same radius of $B(C,1/\eps_0)$. There is a hyperplane $P(C,\omega)\in\mathcal{E}_{\omega_0, A, K, L, \eps_0}$ through $a,b$ (and $C$). For $n>3$ the normal $\omega$, and hence a hyperplane $P(C,\omega)$, are not unique. By conditions of theorem we obtain
$$\|b-C\|\ge1/\eps_0-{\rm dist}(b,K)=1/\eps_0-\delta(K,L)>1/(2\eps_0)>d_L.$$
Hence $L$ does not intersect a line $l(C,\omega)$, and a circular projection $f=f_{C,\,\omega}$
onto $P(C,\omega)$ (see Definition~\ref{D-circ-pro}) is well-defined for $L$.
Since $\|a-C\|=1/\eps_0>d_K$, we conclude that $f$ also is well-defined for $K$.
We have $f(a)=a,\,f(b)=b$ and
\begin{equation*}
  \delta(K,L)=\|a-b\|={\rm dist}(b, f(K))\le\delta(f(L), f(K))
\end{equation*}
that completes a proof of the claim.

In view of (\ref{E-S2}), there is a translation vector $p\,||\,P(C,\omega)$ such that
\begin{equation}\label{E-S17}
 \delta(f(K)+p,\ f(L))\le\eps.
\end{equation}
By Definition~\ref{D-Eomega2}, a circle $S^1_A$ centered at $C$ is orthogonal to $\omega_0$ at a point $A=S_K(\omega_0)$, and $\omega_0$ is a common normal to $K,L$ and $B_K(\omega_0)$ at $A$. Hence $f(A)$ is a boundary point of the images $f(B_K(\omega_0))$ and both $f(K)$ and $f(L)$ (i.e., $f(A)$ belongs to the apparent contours of $B_K(\omega_0)$ and both $K$ and $L$).
Denote by $B'$ a minimal ball supporting $f(B_K(\omega_0))$ at $f(A)$.
By Lemma~\ref{L-Aball2}, the radius of $B'$ satisfies the inequality $R'\le 2R_K(\omega_0)$.
We apply the non-convex version of Lemma~\ref{L-01} for two bodies $M=f(K)+p$ and $N=f(L)$ lying in $P(C,\omega)$ and $O=f(A)$
and obtain that
\begin{equation}\label{E-S16b}
 \delta(f(K), f(L))\le2\sqrt{\eps\,R'}+2\,\eps\le2\sqrt{2\eps\,R_K(\omega_0)}+2\,\eps.
\end{equation}
The desired inequality (\ref{E-S3})(a) is now an immediate consequence of (\ref{E-S16b}) and (\ref{E-S16}). The proof of (b) is founded on Lemma~\ref{L-01b} given below and is similar to~(a).$\,\square$

\begin{prop}\label{P-ab} For any compact subsets $K,L$ of a complete metric space $(M,\rho)$
there is a pair $a\in K,\,b\in L$ such that $\rho(a,b)=\delta(K,L)$, moreover,
either dist$(a,L)=\delta(K,L)$ or dist$(b,K)=\delta(K,L)$.
\end{prop}

 If~$K,L\subset\bbr^n$ are homothetic, that means if $K=\lambda L+p$
 for some $p\in\bbr^n$, $\lambda>0$, we write $K\sim L$.
 The \textit{homothetic deviation} of compact sets $K,L$ is~\cite{gr87},\,\cite{gr97}
 \begin{equation*}
 \delta_h(K,L)=\inf\{\delta(\lambda\,K+p,\,L):\ \lambda>0,\ p\in\bbr^n\}.
 \end{equation*}
Note that if $K\in\mathcal{K}^\eps_i$ then $C+\lambda(K-C)\in\mathcal{K}^{\eps/\lambda}_i$ for any $C\in\bbr^n$ and $\lambda>0$.

\begin{defi}\rm
 If $K,L\subset\bbr^n$ are related each to other by $L{-}C{=}\lambda\,(\phi(K){-}C)$ for some $C\in\bbr^n$, $\phi\in Rot_n(C)$ and $\lambda>0$, we will write $K\overset{C}{\sim} L$.
 The~\textit{homothety-rotational relative} to $C$ \textit{deviation} of compact sets $K,L\subset\bbr^n$ will be called
 \begin{equation*}
 \begin{array}{c}
 \delta_{hC}(K,L)=\inf\{\delta(C+\lambda\,(\phi(K)-C),\,L):\ \phi\in Rot_n(C),\ \lambda>0\}.
 \end{array}
 \end{equation*}
\end{defi}

Since the homotheties $x\to C+\lambda(x-C)$ commute with the rotations $\phi\in Rot_n(C)$, similar estimates of the homothetic deviations $\delta_{h}(K,L)$ or $\delta_{hC}(K,L)$ by means of corresponding deviations $\delta_{h}(f(K), f(L))$ or $\delta_{hC}(f(K), f(L))$ of their circular projections $f=f_{C,\omega}$ can be derived for $K,L\in\mathcal{K}^{\eps}_3$ as well.

\subsection{Auxiliary lemmas}

\begin{lem}[\cite{gr97}]\label{L-01} Let $M,N\subset\bbr^n\ (n\ge2)$ be compact convex bodies and assume that for some $\eps\ge0$
\begin{equation}\label{E-S8}
    \delta(M,N)\le\eps.
\end{equation}
Let $\omega\in S^{n-1}$ be such that $M$ has finite support radius $R_M(\omega)$
in the direction $\omega$.
If $O\in S_N(\omega)$, then the (unique) translate $M'$ of $M$ with the property that $O\in S_{M'}(\omega)$ satisfies the inequality
\begin{equation}\label{E-S9}
    \delta(M',N)\le 2\big(\sqrt{R_M(\omega)}+\sqrt{\eps}\big)\sqrt{\eps}.
\end{equation}
\end{lem}

Lemma~\ref{L-01} and its generalization, Lemma~\ref{L-01b}, are used in the proof of Theorem~\ref{T-01}.
Remark that the convexity of $M,N$ is not used in the proof of Lemma~\ref{L-01}.
We~modify some definitions of \cite{gr97} to convenient for us form.

\begin{defi}\rm Let $K\subset\bbr^n$ be a compact body and $C$ be a point such that $B(C,1/\eps)$ is an outer support ball of $K$. Denote by $r_K(C)$ a minimal real such that a ball $B_K^+(C)=B(C,1/\eps{+}r_K(C))$ contains $K$. We call $S_K^+(C)=\partial B_K^+(C)$ the support sphere of $K$ relative to $C$.
We call $S_K(C)=K\cap S_K^+(C)$ the \textit{support set of $K$ relative to $C$}.
If $S_K(C)$ consists of a point, it is called the \textit{support point of $K$ relative to $C$}. We say that a ball $B\subset\bbr^n$ is an \textit{inner support ball of $K$ relative to $C$} if $K\subseteq B\subset B_K^+(C)$. The~smallest inner support ball of $K$ relative to $C$ is denoted by $B_K(C)$, its radius is denoted by $R_K(C)$. (Clearly if $R_K(C)<1/\eps$, then $S_K(C)$ is a point).
\end{defi}

\begin{lem}\label{L-01b} Let $M,N\subset\bbr^n\ (n\ge2)$ be compact bodies and assume that
\begin{equation}\label{E-S8c}
    \delta(M,N)\le\eps
\end{equation}
for some $\eps\ge0$. Let $C\in\bbr^n$ be such that $B(C,1/\eps_0)$ is an outer support ball of $M$ and the support radius $R_M(C)<1/\eps_0$. If $O\in S_N(C)$, then the (unique) rotation image $M'=\phi(M)$, where $\phi\in Rot_n(C)$, with the property that $\phi(S_M(C))$ is contained in the ray $OC$, satisfies the inequality
\begin{equation}\label{E-S9c}
 \begin{array}{c}
 \delta(M',N)\le \big(2\sqrt{R_M(C)}+3\sqrt{\eps}+\tilde\eps\big)\sqrt{\eps},
 \end{array}
\end{equation}
where
 $\begin{array}{c}
 \tilde\eps=(R_M(C)-\eps)\,\alpha\big(1+\frac{R_M(C)-\eps}{4\,\sqrt{R_M(C)}}\,\alpha\big)
 \end{array}$
 and
 $\begin{array}{c}
  \alpha^2=\frac{\eps_0}{1+\eps_0(r_M(C)-R_M(C))}.
 \end{array}
 $
\end{lem}
\begin{figure*}[ht]
\begin{center}
\begin{tabular*}{\textwidth}{@{\hskip10mm}l@{\extracolsep{\fill}}r@{\hskip15mm}}
 \subfigure[\label{F-Th-4}]{%
\includegraphics[scale= 0.65,angle=0]{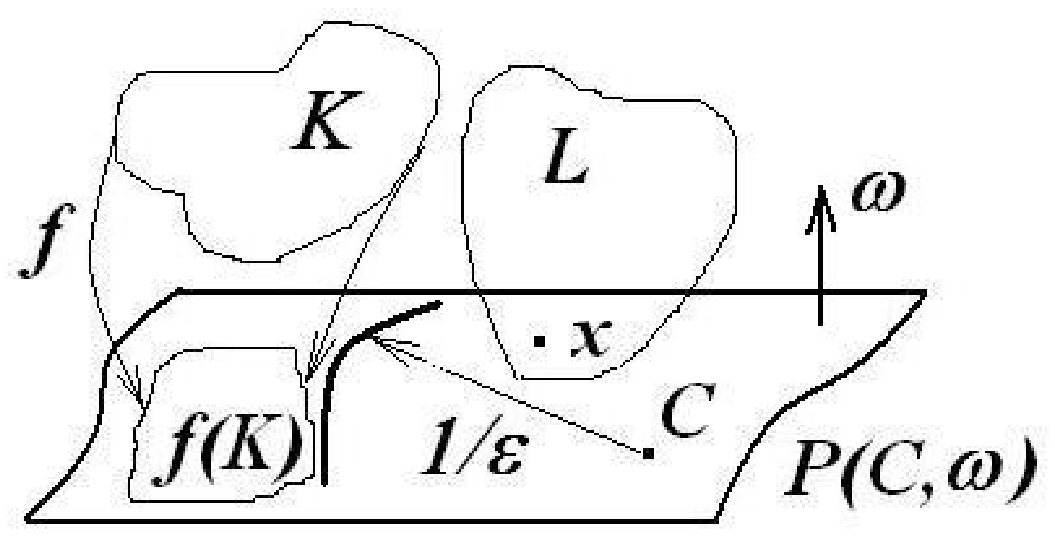}} &
 \subfigure[\label{F-lemma_theta} ]{%
\includegraphics[scale= 0.4,angle=0]{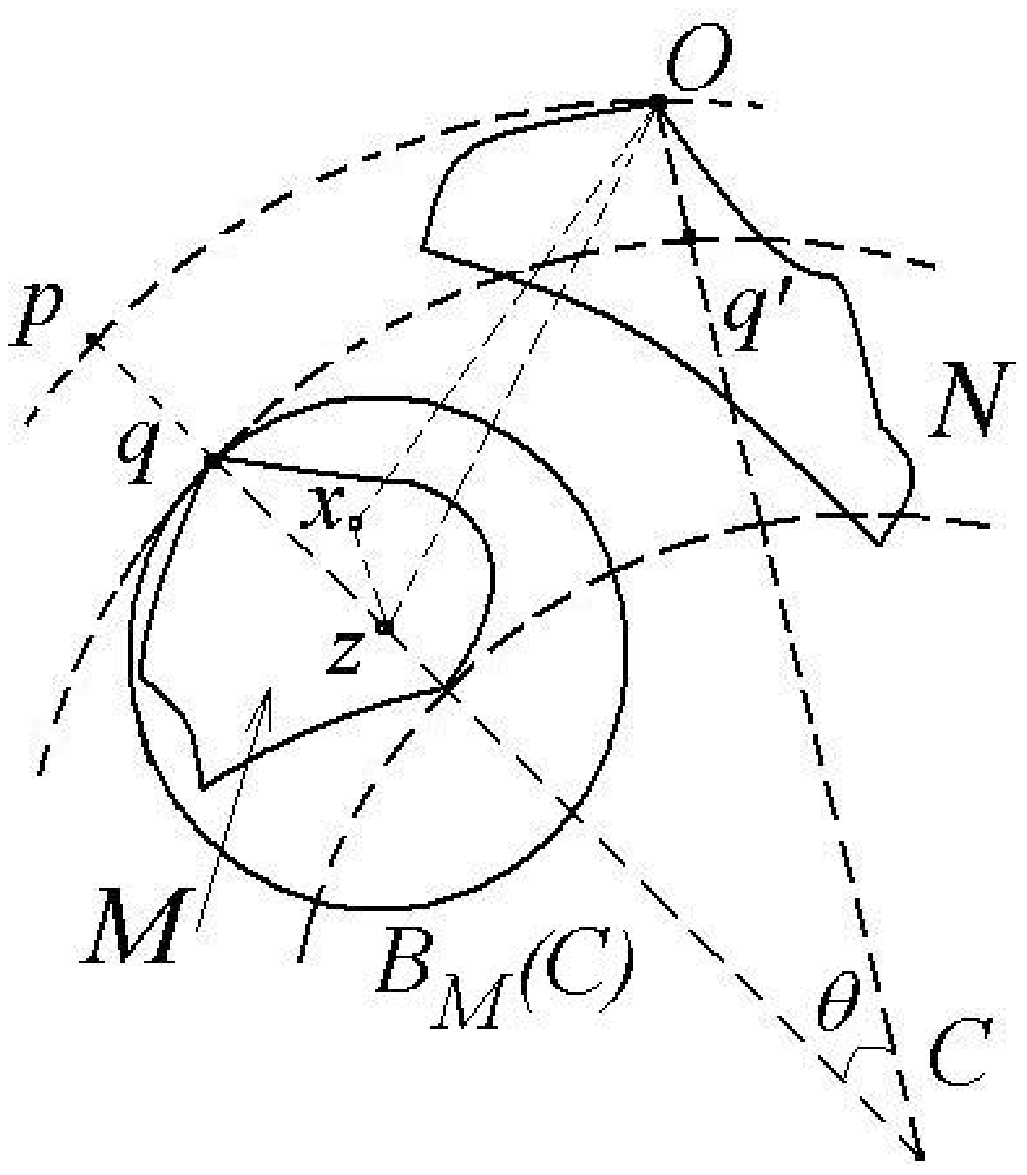}}
\end{tabular*}
\caption{\small (a)~ Proof of Theorem \ref{T-0eps}.\quad
 (b)~Proof of Lemma \ref{L-01b}.}
\end{center}
\end{figure*}
\textbf{Proof}. Denote by $q=S_M(C),\,q'=\phi_C(q)$, Fig~\ref{F-lemma_theta}, and set $\eps'=\|q-q'\|$.
By conditions,  $\|y-\phi_C(y)\|\le\eps'$ holds for all $y\in M$.
Then (\ref{E-S8c}) implies
 $M'\subset N + (\eps'+\eps)\,B(O,1)$ and $N\subset M' + (\eps'+\eps)\,B(O,1)$.
Hence
 $\delta(M',N)\le \eps'+\eps$,
 see (\ref{E-hausd2}), and to prove (\ref{E-S9c}) we only have to show that
\begin{equation}\label{E-S10}
 \eps'\le (2\sqrt{R_M(C)}+2\,\sqrt{\eps}+\tilde\eps)\sqrt{\eps}.
\end{equation}
Let $z$ denotes the center of the ball $B_M(C)$ and let $l$ denotes the line through $z$ and $C$. Clearly, $l$ contains $q$.
Let $p$ denote the intersection point of $l$ and the sphere $S^+_M(C)$
of radius $1/\eps_0+r_M(C)$. Clearly, $\|\phi_C(q)\|=\|q-p\|$.
 As a consequence of (\ref{E-S8c}) and the assumption $O\in N$ there is $x\in M$ such that
\begin{equation}\label{E-g11}
\begin{array}{c}
 \|x\|\le\eps,
\end{array}
 \end{equation}
and, since $\|z-x\|\le R_M(C)$, it follows that
\begin{equation}\label{E-g12}
\begin{array}{c}
 \|z\|\le\|x\|+\|z-x\|\le \eps+R_M(C).
\end{array}
\end{equation}
Next we estimate a "small" angle $\theta=\angle OCq$, see Fig~\ref{F-lemma_theta}. Note that
\begin{eqnarray}\label{E-g12d}
\nonumber
 \|z-C\|=1/\eps_0+r_M(C)-R_M(C),\\
 1/\eps_0+r_M(C)-\eps\le\|C\|=\|p-C\|\le 1/\eps_0+r_M(C)+\eps.
\end{eqnarray}
By cosine theorem for $\triangle COz$ and using (\ref{E-g12}), (\ref{E-g12d}), we obtain
\begin{equation*}
\begin{array}{c}
 \cos\theta=\frac{\|z-C\|^2+\|C\|^2-\|z\|^2}{2\|z-C\|\cdot\|C\|}{\ge}
 \frac{(1/\eps_0{+}r_M(C){-}R_M(C))^2+(1/\eps_0{+}r_M(C){-}\eps)^2-(R_M(C)+\eps)^2}
 {2(1/\eps_0+r_M(C)-R_M(C))(1/\eps_0+r_M(C)+\eps)}\\
 =1-2\,\eps\/\eps_0\frac{1+\eps_0r_M(C)}{[1+\eps_0(r_M(C)-R_M(C))]\,[1+\eps_0(r_M(C)+\eps)]}.
\end{array}
\end{equation*}
Hence
\begin{equation}\label{E-sintheta}
\begin{array}{c}
 \sin^2\frac\theta2=\frac12(1{-}\cos\theta)\le
 \eps\/\frac{\eps_0(1+\eps_0r_M(C))}
 {[1+\eps_0(r_M(C)-R_M(C))]\,[1+\eps_0(r_M(C)+\eps)]}\le\eps\alpha^2.
\end{array}
\end{equation}
Let us assume that
\begin{equation}\label{E-g13}
    \eps\le R_M(C).
\end{equation}
Since $\|p-q\|$ is the distance between the spheres $S^+_M(C)$ and $S(C,\|C\|)$ it follows from (\ref{E-S8c}) that
\begin{equation}\label{E-g14}
    \|p-q\|\le\eps.
\end{equation}
If $p$ is inside $B_M(C)$, then (\ref{E-g14}) implies
\begin{equation*}
    \|z-p\|=R_M(C)-\|p-q\|\ge R_M(C)-\eps,
\end{equation*}
and if $p$ is outside $B_M(C)$, then $\|z-p\|\ge R_M(C)$.
Thus, in either case we have
\begin{equation}\label{E-g15}
    \|z-p\|\ge R_M(C)-\eps.
\end{equation}
Combining (\ref{E-g12}) and (\ref{E-g15}) with the fact that $O,z,p$ are the vertices of an "almost right triangle" we deduce using cosine theorem that
\begin{equation*}
\hskip-1mm\begin{array}{c}
 \|p\|^2{-}2\/\|p\|(R_M(C){-}\eps)\sin\frac\theta2\le
 \|p\|^2{-}2\|p\|{\cdot}\|z{-}p\|\sin\frac\theta2=\|z\|^2{-}\|z{-}p\|^2\\
 \quad \le(\eps+R_M(C))^2-(\eps-R_M(C))^2=4\,R_M(C)\,\eps.
 \end{array}
\end{equation*}
We find roots of the corresponding quadratic equation, and using the inequality $\sqrt{1+x}\le 1+x/2$ and (\ref{E-sintheta}), obtain
\begin{equation*}
\begin{array}{c}
 \|p\|\le(R_M(C)-\eps)\sin\frac\theta2+[\sin^2\frac\theta2(R_M(C)-\eps)^2+4\,R_M(C)\,\eps]^{1/2}\\
 \le(R_M(C){-}\eps)\sin\frac\theta2{+}2\sqrt{R_M(C)\,\eps}
 \big(1{+}\frac{(R_M(C)-\eps)^2}{8\,R_M(C)\,\eps}\sin^2\frac\theta2\big)
 {=}(2\sqrt{R_M(C)}{+}\tilde\eps)\sqrt{\eps}.
 \end{array}
\end{equation*}
This, together with (\ref{E-g14}) yields the conclusion that
\begin{equation*}
 \eps'\le\|p-q\|+\|p\|+\|\phi_C(q)\|\le(2\sqrt{R_M(C)}+2\,\sqrt{\eps}+\tilde\eps)\sqrt{\eps}.
\end{equation*}
If (\ref{E-g13}) is not satisfied, then, using (\ref{E-g11}) and the fact that $q,x\in B_M(C)$, we find
\begin{equation*}
\begin{array}{c}
    \eps'\le \|x\|+\|q-x\|+\|\phi_C(q)\|\le 2\,R_M(C)+2\,\eps
    \le (2\sqrt{R_M(C)}+2\,\sqrt{\eps})\sqrt{\eps}\\
    \le (2\sqrt{R_M(C)}+2\,\sqrt{\eps}+\tilde\eps)\sqrt{\eps}.
 \end{array}
\end{equation*}
Thus, in both cases (\ref{E-S9c}) is valid and this, as already noted, proves the lemma.$\,\square$


\begin{thebibliography}{9999}

\bibitem[1]{rgar06}
 R. Gardner,
{Geometric Tomography}, Cambridge University Press, 2006.

\bibitem[2]{vgol}
 V.\,P. Golubyatnikov,
{Uniqueness Questions in Reconstruction of Multidimensional Objects from Tomography-Type Projection Data}, VSP BV, The Nederlands, 2000.

\bibitem[3]{vgol99}
 V.\,P. Golubyatnikov,
On the unique determination of compact convex sets from their projections. The complex case,
{Siberian Math. J.} 40\,(4), (1999), 678--681.

\bibitem[4]{gr87}
 H. Groemer, Stability theorems for projections of convex sets,
{Israel J. Math.}, 60, No. 2 (1987), 177--190.

\bibitem[5]{gr97}
 H. Groemer, On the determination of convex bodies by translates of their projections,
{Geometriae Dedicata}, 66 (1997), 265--279.

\bibitem[6]{gro}
 M. Gromov, Spaces and questions,
{Geom. Funct. Anal. (GAFA)}, Special Volume, Part I (2000), 118--161.

\bibitem[7]{haef}
 A. Haefliger,
 Quelques remarques sur les applications diff\'{e}rentiables d'une surface dans le plan,
{Ann. Inst. Fourier}, 10 (1960), 47--60.

\bibitem[8]{le80}
 K. Leichtweis,
{Konvexe Mengen}, Berlin, 1980.

\bibitem[9]{mot35}
 T. Motzkin, Sur quelques proprietes characteristiques des ensembles bornes non convexes,
{Rend. Acad. Lincei}, 21 (1935), 773--779.

\bibitem[10]{nic}
 C.\,P. Niculescu, L.-E. Persson,
{Convex Functions and Their Applications: a Contemporary Approach},
Canad. Math. Series Books in Mathematics, Springer, 2006.

\bibitem[11]{resh4}
Yu.\,G. Reshetnyak (Ed.),
{Geometry IV, Non-regular Riemannian Geometry. Series: Encyclopaedia of Mathematical Sciences, 70}, Springer, 1993.

\bibitem[12]{resh56}
 Yu.\,G. Reshetnyak,
 {On one generalization of convex surfaces}. Math. Sbornik, 1956, 40 (3), 381--398 (Russian).

\bibitem[13]{Rov_2006}
 V.\,Yu. Rovenski, On $(k,\varepsilon)$-saddle submanifolds of Riemannian manifolds,
{Geometriae Dedicata}, 121 (2006), 187--203.

\bibitem[14]{Shar_94}
V.\,A. Sharafutdinov,
{Integral Geometry of Tensor Fields}, VSP BV, Utrecht, The Nederlands, 1994.

\bibitem[15]{toj_06}
R. Tojeiro, Riemannian $G$-manifolds as Euclidean submanifolds,
{Revista de la Uni\'{o}n Matem\'{a}tica Argentina}, 47, No. 1 (2006), 73-–83.

\bibitem[16]{top}
V.\,A. Toponogov,
{Differential Geometry of Curves and Surfaces. A Concise Guide}. With the editorial assistance of V.Rovenski, Birkh\"{a}user,~2006.

\end{thebibliography}
\end{document}